\documentclass[11pt]{article}

\usepackage{amsmath}
\usepackage{amsfonts}
\usepackage{pgfplots}
\usepackage{amssymb}
\usepackage{tikz}
\usetikzlibrary{positioning}
\usetikzlibrary{matrix,arrows,decorations.pathmorphing}
\usetikzlibrary{patterns}
\usetikzlibrary{through,calc,3d}
\usetikzlibrary{plotmarks}
\usepackage[titletoc]{appendix}
\usepackage[T1]{fontenc}
\usepackage[numbers,square,sort&compress]{natbib}
\usepackage{subfigure}                  
 \usepackage[affil-it]{authblk}

\newtheorem{theorem}{\bf Theorem}[section]

\newcommand{\pfrac}[2]{\ensuremath{\dfrac{\partial #1}{\partial #2}}}
\newcommand{\inner}[2]{\big< #1 ,  #2 \big>}
\newcommand{\bm}[1]{\ensuremath{\mathbf{#1}}}
\newcommand{\bs}[1]{\ensuremath{\boldsymbol{#1}}}

\newcommand{\udot}[1]{\dot{\bm{u}}_{#1}}
\newcommand{\ydot}[1]{\dot{\bm{y}}_{#1}}
\renewcommand{\u}[1]{\bm{u}_{#1}}
\newcommand{\y}[1]{\bm{y}_{#1}}

\DeclareMathOperator{\tr}{trace}

\usepackage[font=scriptsize]{caption}

\begin{document}

\title{Computing Sensitivities in Evolutionary Systems: A Real-Time Reduced Order Modeling Strategy}
\author[1]{Michael Donello}
\author[2]{Mark Carpenter}
\author[1]{Hessam Babaee\thanks{Corresponding author. Email:h.babaee@pitt.edu.}}
\affil[1]{\footnotesize Department of Mechanical Engineering and Materials Science, University of Pittsburgh}
\affil[2]{\footnotesize NASA Langley Research Center, Hampton, VA}
\date{}
\maketitle

\begin{abstract}
    We present a new methodology for computing sensitivities in evolutionary systems using a model driven low-rank approximation. To this end, we formulate a variational principle that seeks to minimize the distance between the time derivative of the reduced approximation and sensitivity dynamics. The first order optimality condition of the variational principle leads to a system of closed-form evolution equations for an orthonormal basis and corresponding sensitivity coefficients. This  approach allows for the computation of sensitivities with respect to a large  number of  parameters in an accurate and tractable manner by extracting correlations between different sensitivities on the fly. The presented method requires solving  \emph{forward} evolution equations,  sidestepping  the restrictions imposed  by \emph{forward/backward} workflow of adjoint  sensitivities. For example, the presented method, unlike adjoint equation, does not impose any I/O load  and can be used in applications in which real time sensitivities are of interest.  We demonstrate the utility of the method for three test cases: (1) computing sensitivity with respect to model parameters in the R\"{o}ssler system (2) computing sensitivity with respect to an infinite dimensional forcing parameter in the chaotic Kuramoto-Sivashinsky equation and (3) computing sensitivity with respect to reaction parameters for species transport in a turbulent reacting flow. 
\end{abstract}


\section{Introduction}

Sensitivity analysis is required in a diverse set of evolutionary systems that are governed by differential equations in the form of $\dot{\bm{v}} = \bm{g}(\bm{v};\bs \alpha)$, where $\dot{(\sim)}=d(\sim)/dt$, $\bm{v} \in \mathbb{R}^n$ is the state space variable and $\bs \alpha \in \mathbb{R}^d$ is the design space. These sensitivities, denoted by $\bm{v}'_i = \partial \bm{v}/\partial \alpha_i$, $i=1, \dots, d$,  are needed in numerous applications such as gradient-based optimization \cite{G00,G99}, optimal control \cite{BMT01}, grid adaptivity \cite{DK19}, and parameter identification \cite{EC11}, to name a few. The sensitivities are commonly computed via finite difference (FD), or by directly solving a sensitivity equation (SE) or adjoint equation (AE).  The computational cost of using FD or SE scales linearly with the number of  parameters -- making them impracticable when sensitivities with respect to a large number of parameters are needed.  On the other hand, the computational cost of solving AE is independent of the number of parameters as it requires solving a single ordinary/partial differential equation (ODE/PDE).

 While AE is certainly a preferred approach for computing sensitivities for stationary problems, for time-dependent problems, the  \emph{forward-backward} workflow of adjoint solver can pose several challenges. In particular,  solving AE imposes a significant storage cost as the AE must be solved backward in time. On the other hand,   the    adjoint operator utilizes the  forward time-resolved solution of the nonlinear dynamical system, i.e. $\bm{v}$. As a result, the dynamical system must be solved forward in time and its time-resolved solution must be stored. The adjoint solver is then solved backward in time, in which the nonlinear state is read from the disk at every time step. This workflow is not adequate for problems where real-time sensitivities are required, e.g. grid adaptivity for time-dependent problems \cite{DK19}. Moreover,  for  high dimensional dynamical systems, i.e. $n \sim \mathcal{O}(10^{10})$, the imposed  I/O operations in the AE workflow could lead to  insurmountable limitations.

The   I/O limitations continue to become more restrictive in the future high performance computing architectures and it is one of the major challenges in the transition from current sub-petascale and petascale computing to exascale computing \cite{AM_DOE_14}.  For example,  in high fidelity simulations of turbulent reactive flows, the solution can only be stored at every 400th time step in order to maintain I/O overhead at a reasonable level, while important events such as the ignition kernel
 occur rapidly on the order of 10 simulation time steps \cite{BAB12}. Storing the time-resolved solution for these problems is required for AE and it is currently exceeding the acceptable  I/O levels and this trend continues to become even more unfavorable for exascale computing.   This alone
 gives rise to a growing need for algorithms that can accurately compute sensitivities
 while minimizing or eliminating I/O requirements, and this is one of the motivations of the method presented in this paper.   \\

Sensitivities of a dynamical system with respect to different parameters are often highly correlated and therefore they are amenable to low-rank approximations. To this end,   a new  low-dimensional model was recently presented in \cite{Babaee_PRSA} that can  describe transient instabilities in high-dimensional nonlinear dynamical systems. This approach is based on a time-dependent basis known as the optimally time-dependent (OTD) modes. The evolution equation for the OTD modes is obtained by minimizing the functional
 \begin{equation}
\mathcal{F}(\dot {\bm{u}}_1,\dot{\bm{u}}_2, \dots,\dot{\bm{u}}_r) 
 =\sum_{i=1}^{r}\big\Vert \dot {\bm{u}}_i-\bm{L}({\bm{v}(t),t)\bm{u}_i(t)}\big\Vert ^{2},\label{eq:functional}
\end{equation}
subject to the orthonormality of the OTD modes, i.e. $\bm{u}_i^T\bm{u}_j = \delta_{ij}$, where $\bm{u}_i \in \mathbb{R}^n, i=1, \dots, r$ are the OTD modes.  In the above functional, $\left\Vert \bm{u}\right\Vert ^{2}= \bm{u}^T\bm{u}$ and $\bm{L}(\bm{v}(t),t):=\nabla_{\bm{v}}\bm{g}$ is the instantaneous linearized operator. The optimality condition of the above variational principle leads to a closed form evolution equation for the OTD subspace:  $\dot{\bm{U}} = (\bm{I}-\bm{U}\bm{U}^T)\bm{L} \bm{U}$, where  $\bm{U}  = [\bm{u}_1 | \bm{u}_2 | \dots  |\bm{u}_r ] \in \mathbb{R}^{n\times r}$ and $\bm{I}\in \mathbb{R}^{n \times n}$ is the identity matrix.   It was shown later that the OTD subspace converges exponentially fast  to the eigendirections of the Cauchy–Green tensor associated with the most intense finite-time instabilities \cite{BFHS17}. In this sense, the OTD reduction can be interpreted as a low-rank subspace that approximates the evolution of perturbed initial condition in all directions of the phase space. One of the computational advantages of OTD is that it only requires solving forward equations. Moreover, the computational complexity of solving OTD reduction scales linearly with respect to the number of modes. The OTD  has also been used for flow control \cite{BMS18}, building precursors for bursting phenomena  \cite{FS16}   as well as detection of  edge manifolds in infinite-dimensional dynamical systems \cite{BDSH20}.  We also note that the time-dependent bases have also been developed in the context of stochastic reduced order modeling; see for example \cite{SL09,CHZI13,MNZ15,Babaee:2017aa,PB20}. 

Our objective in this paper is to approximate sensitivities with respect to a large number of parameters using  forward low-rank  systems similar to OTD.  In particular we seek to reduce the computational cost of solving forward sensitivity equations by exploiting the \emph{correlations} between various sensitivities.    However,  OTD is not adequate when   applied to systems subject to perturbations in a parametric space.  These perturbations are governed by the forced linear  sensitivity equation:  $\partial \bm{v}'_i/\partial t = \bm{L} \bm{v}'_i + \partial \bm{g}/\partial \alpha_i $, and in general, the OTD subspace is not an optimal basis for the evolution of $\bm{v}'_i$. To this end, we present a new approach based on time-dependent basis for solving  time-varying linear systems forced by a high-dimensional function. 



The contributions of this paper are twofold: (i) we present a new variational principle, whose optimality conditions lead to forward real-time  low rank evolution equations for the approximation of the forced sensitivity equation.  We coin this approach ``forced OTD", which we will simply refer to as f-OTD. (ii) We extend the application of the presented method to compute \emph{tensor-like} sensitivities. An example of tensor-like sensitivities is in reactive flows where the goal is to compute the sensitivity of $n_s$ species with respect to $m$ parameters. In these systems, the full sensitivities can be represented as a third-order tensor  where the first dimension is the number of grid points, the second dimension represents species ($n_s$), and the third dimension represents the parameters ($m$). We show that with a single set of orthonormal modes, we can approximate  sensitivities by exploiting correlations between \emph{all}  sensitivities. We compare the computational cost of f-OTD against adjoint based sensitivities where one adjoint variable for each species must be solved \cite{BOR15,LCR19,LLC20}.    
We show how the presented approach can be used for computing sensitivities with respect to a large number of parameters by  solving forward low-rank evolution equations without the need to store the state variables.

In the sections that follow, we present the formulation of the f-OTD method and demonstrate a number of  outcomes. We start in Section 2 with the variational principle  whose optimality conditions lead to a set of closed form evolution equations for a low rank approximation of the forced sensitivity equation. In Section 3, we present three demonstration cases: (1) sensitivity with respect to model parameters in the R\"{o}ssler system (2) sensitivity with respect to an infinite dimensional forcing parameter in the chaotic Kuramoto-Sivashinsky equation and (3) sensitivity with respect to reaction parameters for species transport in a turbulent reacting flow. In Section 4, we present the conclusions and implications of our work. 





\section{Methodology}\label{sec:Methodology}
\subsection{Preliminaries}
We denote $\bm{u}(x,t)$ to be a time dependent field variable. We denote the spatial domain as $D \subset \mathbb{R}^m$, where $m=$ 1, 2, or 3. The spatial coordinate is denoted by $x\in D$ and the function is evaluated at time $t$. We define the inner product of  functions $\bm{u}(x,t)$ and $\bm{v}(x,t)$ as 
\begin{equation*}
    \inner{\bm{u}(x,t)}{\bm{v}(x,t)} = \displaystyle \int\limits_{D} \bm{u}(x,t) \bm{v}(x,t)  dx
\end{equation*}
and  the $L_2$ norm  induced by the above inner product  is given as
\begin{equation*}
    \left\Vert\bm{u}(x,t)\right\Vert_{2} = \inner{\bm{u}(x,t)}{\bm{u}(x,t)}^{\frac{1}{2}}.
\end{equation*}
We introduce a quasimatrix notation to represent a set of functions in matrix form, and denote the quasimatrix $\bm{U}(x,t)\in \mathbb{R}^{\infty\times r}$  as \cite{BT04}:
\begin{equation*}
    \bm{U}(x,t) = \bigg[\bm{u}_1(x,t)\  \Big| \ \bm{u}_2(x,t) \ \Big| \  \dots \ \Big| \ \bm{u}_d(x,t)   \bigg]_{\infty \times r},
\end{equation*}
where the first dimension is infinite and represents the continuous state space contained by $D$ and the second dimension is discrete.  Similarly, we use the term quasitensor for tensors whose first dimension is infinity. For example, $\bs{\mathcal{T}} \in \mathbb{R}^{\infty \times r_1 \times r_2}$ is a third-order quasitensor.  Following this definition, we define the inner product between quasimatrices $\bm{U}(x,t)\in\mathbb{R}^{\infty\times r}$ and $\bm{V}(x,t)\in \mathbb{R}^{\infty\times d}$  as 
\begin{equation*}
    \bm{S}(t) = \inner{\bm{U}(x,t)}{\bm{V}(x,t)},
\end{equation*}
where $\bm{S}(t)\in\mathbb{R}^{r\times d}$ is a matrix with components $S_{ij}(t) = \inner{\bm{u}_i(x,t)}{\bm{v}_j(x,t)}$, where $\bm{v}_j(x,t)$ is the $j$th column of $\bm{V}(x,t)$. The discrete analogue of this operation is the matrix multiplication, $\bm{U}(t)^T \bm{V}(t)$, where $\bm{U}(t)\in\mathbb{R}^{n\times r}$ and $\bm{V}(t)\in\mathbb{R}^{n\times d}$ are space discrete with $n$ grid points. Correspondingly,  the Frobenius norm of a quasimatrix is defined as: 
\begin{equation*}
    \Big \Vert\bm{U}(x,t)\Big \Vert_{F} = \sqrt{\tr \inner{\bm{U}(x,t)}{\bm{U}(x,t)}}.
\end{equation*}
Similarly, we define the inner product between quasimatrix $\bm{U}(x,t)$ and function $\bm{v}(x,t)$ as:
\begin{equation*}
    \bm{g}(t) = \inner{\bm{U}(x,t) }{\bm{v}(x,t)},
\end{equation*}
where $\bm{g}(t)=(g_1(t),g_2(t),\dots,g_r(t))^T \in \mathbb{R}^{r\times 1}$ is a vector with components $g_i(t)=\inner{\bm{u}_i(x,t)}{\bm{v}(x,t)}$. The discrete analogue of this operation is the matrix vector multiplication, $\bm{U}(t)^T \bm{v}(t)$, where $\bm{v}(t)\in\mathbb{R}^{n\times 1}$ is space discrete with $n$ grid points. Finally, we define multiplication between a quasimatrix and a vector
\begin{equation*}
   \bm{c}(x, t) =  \bm{U}(x,t)\bm{b}(t),
\end{equation*}
where $\bm{b}(t) = (b_1(t),b_2(t),\dots,b_r(t))^T\in\mathbb{R}^{r\times 1}$ is an arbitrary vector and $\bm{c}(x,t)\in\mathbb{R}^{\infty\times 1}$ is a function given by $\bm{c}(x, t)=b_i(t)\bm{u}_i(x,t)$. We use index notation and the same indexes imply summation. 

We consider the nonlinear partial differential equation (PDE) for the evolution of $\bm{v}(x,t)$:
\begin{equation}\label{eq:gennon}
    \pfrac{\bm{v}(x,t)}{t} = \mathcal{N}\left((\bm{v}(x,t); \bs \alpha \right), \quad t\in [0, T_f]
\end{equation}
where $\mathcal{N}$ is in general a nonlinear differential operator. Our goal is to compute the sensitivity of $\bm{v}(x,t)$ with respect to the design parameters $\bs{\alpha}$, which can either be infinite-dimensional, i.e., a function $\bs{\alpha} = \bs{\alpha}(x,t)$, or finite-dimensional, i.e., a vector $\bs{\alpha} = (\alpha_1, \alpha_2, \dots, \alpha_d)$. For the sake of simplicity in the exposition we consider the finite-dimensional parametric space. Differentiating Equation \ref{eq:gennon} with respect to design parameter $\alpha_i$ leads to an evolution equation for the sensitivity of the dynamical system:
\begin{equation}\label{eqn:sensitivity}
    \pfrac{\bm{v}_{i}'(x,t)}{t} = \mathcal{L}\left(\bm{v}_{i}'(x,t)\right) + \bm{f}'_{i}(x,t;\bs{\alpha}),
\end{equation}
where $\bm{v}_{i}' = \partial \bm{v}/\partial \alpha_i$ is the sensitivity of $\bm{v}(x,t)$ with respect to $\alpha_i$, $\mathcal{L} (\sim ) = \partial \mathcal{N}/ \partial \bm{v} (\sim )$ is the linearized operator, and $\bm{f}'_{i} = \partial\mathcal{N}/ \partial \alpha_i$ is the forcing term. 

\subsection{Variational Principle for Reduced Order Modeling} 
\begin{figure}[tbp]
    \centering
    \includegraphics[width=1\textwidth]{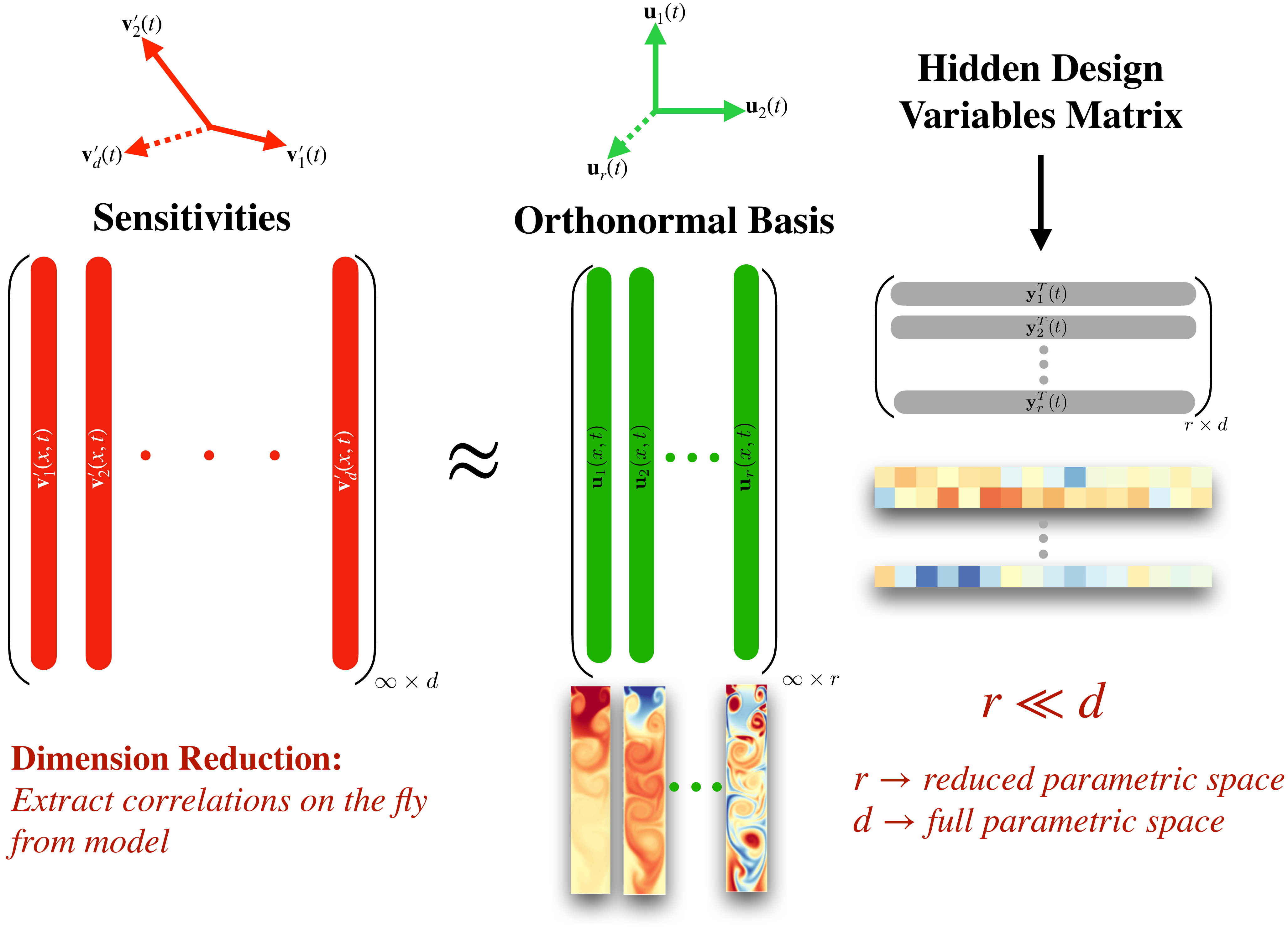}
    \caption{Overview of the reduced order modeling strategy. Shown on left in red is the full dimensional system of sensitivities that we seek to model using the f-OTD low-rank approximation. Shown on right is the low-rank approximation which consists of a set of temporally evolving orthonormal modes (green) and hidden design variables (gray). The hidden design variables are coefficients that map the orthonormal basis to each sensitivity in the full-dimensional system. That is, each of the $d$ sensitivities are approximated as a linear combination of the $r$ orthonormal modes, where $r \ll d$. It is important to note that the orthonormal basis and hidden design variables are model-driven and evolve based on the linear sensitivity dynamics. Thus, the proposed method only requires solving a system of $r$ PDEs and $r$ ODEs   for the modes and coefficients, respectively.}
    \label{fig:schematic}
\end{figure}

Different  sensitivities in a dynamical system tend to be highly correlated at any given time and therefore, these sensitivities can potentially be approximated effectively by a low rank time-dependent subspace.  In this section, we present a real-time reduced order modeling strategy  that aims to extract this subspace and utilize it for building sensitivity ROMs. In particular, we present a variational principle, whose first-order optimality conditions lead to the evolution equations  of a time-dependent subspace and its coefficients. We estimate the sensitivities  using the low-rank decomposition:
\begin{equation}
   \bm{V}'(x,t) =  \bm{U}(x,t)\bm{Y}(t)^T + \bm{E}(x,t),
\end{equation}
where $\bm{V}'(x,t) = \big[\bm{v}'_1(x,t)\  \big| \ \bm{v}'_2(x,t) \ \big| \  \dots \ \big| \ \bm{v}'_d(x,t)   \big]_{\infty \times d}$ is the sensitivities quasimatrix, $\bm{U}(x,t) = \big[\bm{u}_1(x,t)\  \big| \ \bm{u}_2(x,t) \ \big| \  \dots \ \big| \ \bm{u}_r(x,t)   \big]_{\infty \times r}$ is  a quasimatrix representing a rank-$r$ time-dependent orthonormal basis, in which $\inner{\bm{u}_i(x,t)}{\bm{u}_j(x,t)} = \delta_{ij}$, $\bm{Y}(t)=\big[\bm{y}_1(t)\  \big| \ \bm{y}_2(t) \ \big| \  \dots \ \big| \ \bm{y}_r(t)   \big]_{d\times r}$ is the coefficient matrix, and $\bm{E}(x,t) \in \mathbb{R}^{\infty \times d}$ is the approximation error.  The f-OTD decomposition is shown schematically in Figure \ref{fig:schematic}. 

We formulate a variational principle with control parameters $\dot{\bm{U}}(x,t)$ and $\dot{\bm{Y}}(t)$, that seeks to optimally update the subspace $\bm{U}(x,t)$ and its coefficients $\bm{Y}(t)$ by minimizing the distance between the time derivative of the low-rank approximation and the full-dimensional  sensitivity dynamics:
\begin{equation}\label{eq:functional}
\mathcal{F}(\dot{\bm{U}}(x,t), \dot{\bm{Y}}(t)) = \left\Vert \pfrac{(\bm{U}(x,t) \bm{Y}(t)^T)}{t}  - \mathcal{L}\left(\bm{U}(x,t) \bm{Y}(t)^T\right)-\bm{F}'(x,t;\bs{\alpha}) \right\Vert_{F}^{2}.
\end{equation}
where $\bm{F}'(x,t) = \big[\bm{f}'_1(x,t)\  \big| \ \bm{f}'_2(x,t) \ \big| \  \dots \ \big| \ \bm{f}'_d(x,t)   \big]_{\infty \times d}$.
Taking the time derivative of the orthonormality condition leads to the following constraint for the minimization problem: 
\begin{equation}\label{eq:orthodot}
\inner{\dot{\bm{u}}_i(x,t)}{\bm{u}_j(x,t)} + \inner{\bm{u}_i(x,t)}{\dot{\bm{u}}_j(x,t)} = 0.
\end{equation}
We denote $\bs{\phi}_{ij}(t) = \inner{\bm{u}_i(x,t)}{\dot{\bm{u}}_j(x,t)}$, in which $\bs{\Phi}(t)=[\phi_{ij}(t)]\in \mathbb{R}^{r\times r}$. It is easy to see that $\bs{\Phi}(t)$ must be a skew-symmetric matrix in order to satisfy Equation \ref{eq:orthodot}, i.e., $\bs{\phi}_{ji}(t) = -\bs{\phi}_{ij}(t)$. Incorporating this constraint leads to the following unconstrained optimization problem functional:
\begin{align}\label{Eqn:min_princ}
 \mathcal{G}(\dot{\bm{U}}(x,t), \dot{\bm{Y}}(t),\lambda(t)) &=\left\Vert  \pfrac{(\bm{U}(x,t) \bm{Y}(t)^T)}{t}  - \mathcal{L}\left(\bm{U}(x,t)\right) \bm{Y}(t)^T-\bm{F}'(x,t;\bs{\alpha})    \right\Vert_F^{2}\\ \nonumber
 &+ \sum_{i,j=1}^r \lambda_{ij}(t)  \big( \inner{\bm{u}_i(x,t)}{\dot{\bm{u}}_j(x,t)} - \bs{\phi}_{ij}(t) \big),
 \end{align}
 where $\lambda(t) = [\lambda_{ij}(t)] \in \mathbb{R}^{r \times r}$ are Lagrange multipliers. To derive the optimality conditions, we follow a procedure similar to the one that was recently presented in \cite{B19}.  In Appendix \ref{app:A}, we show that minimizing the above functional with respect to $\dot{\bm{U}}(x,t)$ and $\dot{\bm{Y}}(t)$ leads to closed form evolution equations for the modes and corresponding sensitivity coefficients:
 \begin{align}\label{eqn:evol_modes}
     \pfrac{\bm{u}_i(x,t)}{t} &= \mathcal{L}\left(\bm{u}_i\right) - \inner{\bm{u}_j}{\mathcal{L}\left(\bm{u}_i\right)}\bm{u}_j+ \big[ \bm{F}'\bm{y}_k  - \inner{\bm{u}_j}{\bm{F}'\bm{y}_k}\bm{u}_j\big] C_{ik}^{-1} - \phi_{ij}\u{j},
\end{align}
\begin{align}\label{eqn:evol_coeff}
     \frac{d\bm{y}_i(t)}{dt} &= \inner{\bm{u}_i}{\mathcal{L}(\bm{u}_j)} \bm{y}_j + \inner{\bm{F}'}{\bm{u}_i} - \phi_{ij}\y{j},
 \end{align}
where $\bm{C}(t)=[C_{ik}(t)] \in \mathbb{R}^{r\times r}$ is the low-rank \emph{correlation} matrix, in which
$C_{ik}(t) = \bm{y}_i(t)^T \bm{y}_k(t)$. These equations are initialized by solving Equation \ref{eqn:sensitivity} for a single time step and computing the singular value decomposition (SVD) of $\bm{V}'(x,t=\Delta t)$, such that $\bm{U}(x,t=\Delta t)$ contains the first $r$ left singular vectors and $\bm{Y}(t=\Delta t)$ is the matrix multiplication of the first $r$ right singular vectors and singular values; see Section \ref{sec:approx_error}. We show in Section \ref{Sec:Eq} that the skew symmetric matrix $\phi_{ij}$ can be taken to be zero, i.e., $\phi_{ij}=0$.

   In the following, we make several observations about Equations \ref{eqn:evol_modes} and \ref{eqn:evol_coeff}: (i) Equation \ref{eqn:evol_modes} determines the evolution of the f-OTD subspace. For $\phi_{ij}=0$, the right hand side of Equation \ref{eqn:evol_modes} is equal to  the projection of $\mathcal{L}\left(\bm{U}\right) + \bm{F}\bm{Y}\bm{C}^{-1}$ onto the complement of the space spanned by $\bm{U}$.  Therefore, if  $\mathcal{L}\left(\bm{U}\right) + \bm{F}\bm{Y}\bm{C}^{-1}$ is in the span of $\bm{U}$, the f-OTD subspace does not evolve, i.e., $\dot{\bm{U}}=\bm{0}$. However, when  $\mathcal{L}\left(\bm{U}\right) + \bm{F}\bm{Y}\bm{C}^{-1}$ is not in the span of $\bm{U}$, the f-OTD subspace evolves optimally to follow the right hand side.  Equation \ref{eqn:evol_coeff} is the f-OTD reduced order model (ROM) that determines the evolution of the sensitivities within the f-OTD subspace. 
   (ii) We observe that if we set $\bm{F}'(x,t)=\bm{0}$ in the above equations, we recover the  OTD evolution equations presented in \cite{Babaee_PRSA}. However, unlike the OTD equations, where the evolution of the OTD modes are independent of the evolution of the coefficients  ($\bm{Y}$), there is a two-way nonlinear coupling between the f-OTD evolution equations for $\bm{U}$ and $\bm{Y}$. (iii) From the above equations, it is clear to see that f-OTD extracts the low-rank approximation directly from the sensitivity evolution equation. In that sense, it is different from data-driven low-rank approximations such as proper orthogonal decomposition \cite{A91,ABGA15,PTB20} or dynamic mode decomposition \cite{S10,LKB18}, in which the low rank subspace is extracted from preexisting data. The need to generate data simply does not exist in the f-OTD workflow. (iv) The computational cost of solving the f-OTD Equations \ref{eqn:evol_modes} and \ref{eqn:evol_coeff} is roughly equivalent to that of solving $r$ forward sensitivity equations. This is because the evolution of the f-OTD modes described by Equation \ref{eqn:evol_modes} inherits the same differential operators from the sensitivity equation and the computational cost of solving each f-OTD mode is roughly equivalent to that of solving the sensitivity equation for a single parameter. Equation \ref{eqn:evol_coeff} is an ODE and therefore its computational cost is negligible compared to the f-OTD modes, which are governed by a PDE. 

\subsection{Equivalence}\label{Sec:Eq}
It is important to note that the choice of $\phi_{ij}$ in Equations \ref{eqn:evol_modes} and \ref{eqn:evol_coeff} is not unique, and any skew-symmetric matrix yields an equivalent reduction.
Similar to the OTD equations \cite{Babaee_PRSA}, we choose $\phi_{ij}=0$, which corresponds to the dynamically orthogonal (DO) condition.  This property is summarized in the theorem below.
\begin{theorem} 
Let $\{\bm{U}(x,t),\bm{Y}(t)\}$ and $\{\tilde{\bm{U}}(x,t),\tilde{\bm{Y}}(t)\}$ represent two reductions that satisfy Equations \ref{eqn:evol_modes} and \ref{eqn:evol_coeff} with corresponding skew-symmetric matrices $\bs{\Phi}(t)$ and $\tilde{\bs{\Phi}}(t)$, respectively. If the reductions are equivalent at $t=0$, i.e. they are initially related by an orthogonal rotation matrix $\bm{R}_0\in\mathbb{R}^{r\times r}$ as $\bm{U}(x,0)=\tilde{\bm{U}}(x,0)\bm{R}_0$ and $\bm{Y}(0)=\tilde{\bm{Y}}(0)\bm{R}_0$, then the two reductions will remain equivalent for $t>0$ with rotation matrix $\bm{R}(t)$ governed by $\dot{\bm{R}}=\bm{R}\bs{\Phi}-\tilde{\bs{\Phi}}\bm{R}$. 
\end{theorem}
For proof of the above theorem see Appendix \ref{app:B}.

\subsection{Exactness of f-OTD}
For the case where the full sensitivity quasimatrix is of rank $d$, the rank $d$ f-OTD equations are exact. To show this, we start by considering an arbitrary perturbation subspace, $\bm{V}'(x,t)\in\mathbb{R}^{\infty\times d}$, governed by the quasimatrix form of Equation \ref{eqn:sensitivity}:
\begin{align*}
    \pfrac{\bm{V}'}{t} = \mathcal{L}(\bm{V}') + \bm{F}'(x,t), \quad \bm{V}'(x,0) = \bm{V}'_0(x),
\end{align*}
 where columns of $\bm{V}'(x,t)$ are independent, i.e. $\inner{\bm{v}_i'}{\bm{v}_j'}=0$ if $i\neq j$,
and the evolution of an orthonormal subspace, $\bm{U}(x,t)\in\mathbb{R}^{\infty\times d}$, governed by the quasimatrix form of Equation \ref{eqn:evol_modes}:
\begin{align*}
    \pfrac{\bm{U}}{t} = \mathcal{L}(\bm{U}) - \bm{U}\bm{L}_r(t) + \left( \bm{F}'\bm{Y} - \bm{U}\inner{\bm{U}}{\bm{F}'\bm{Y}} \right)\bm{C}^{-1}, \quad \bm{U}(x,0) = \bm{U}_0(x).
\end{align*}
The corresponding matrix of sensitivity coefficients are governed by the matrix form of Equation \ref{eqn:evol_coeff} as:
\begin{align*}
    \frac{ d\bm{Y}}{dt} = \bm{Y}\bm{L}_r^T + \inner{\bm{F}'}{\bm{U}}, \quad \bm{Y}(0)=\bm{Y}_0,
\end{align*}
where $\bm{L}_r(t) = \inner{\bm{U}(x,t)}{\mathcal{L}(\bm{U}(x,t))}$ is the $r\times r$ low-rank linear operator. We can show that if the two subspaces are initially equivalent, i.e., $\bm{U}_0(x)$ can be mapped to $\bm{V}'_0(x)$ via the linear transformation $\bm{Y}_0^T$, then $\bm{V}'(x,t)$ and $\bm{U}(x,t)$ remain equivalent for all time $t$ and are related by the linear transformation $\bm{Y}(t)^T$. This leads to the following theorem:

\begin{theorem}
Let $\bm{V}'(x,t) \in \mathbb{R}^{\infty\times d}$ be an arbitrary subspace evolved by the linear dynamics of Equation \ref{eqn:sensitivity}, and $\bm{U}(x,t)\in \mathbb{R}^{\infty\times d}$ be an orthonormal subspace evolved by Equation \ref{eqn:evol_modes}. If initially $\bm{V}'_0(x)$ and $\bm{U}_0(x)$ are equivalent, i.e. $\bm{V}'_0(X)=\bm{U}_0(x)\bm{Y}_0^T$, then the perturbation subspace can be exactly determined via the linear transformation $\bm{V}'(x,t)=\bm{U}(x,t)\bm{Y}(t)^T$ for all time $t$, where $\bm{Y}(t)$ is governed by Equation \ref{eqn:evol_coeff}. 
\end{theorem}
For a detailed proof of the theorem, refer to Appendix \ref{app:C}.

\subsection{Approximation error}\label{sec:approx_error}
The approximation error of estimating sensitivities using f-OTD can be expressed as $e(t) = \| \bm{V}'(x,t) - \bm{U}(x,t)\bm{Y}(t)^T \|_F$. This error can be properly analyzed and  better understood by considering two types of error:  (i) the resolved  error, denoted by $e_r(t)$ and (ii) the unresolved error, denoted by $e_u(t)$. The resolved error is the discrepancy between approximating the sensitivities with rank-$r$ f-OTD and the optimal rank-$r$ approximation: $e_r(t) = \|  \bm{U}(x,t)\bm{Y}(t)^T - \tilde{\bm{U}}(x,t)\tilde{\bm{Y}}(t)^T \|_F$, where $\tilde{\bm{U}}(x,t) \in \mathbb{R}^{\infty \times r}$  and $\tilde{\bm{Y}}(t) \in \mathbb{R}^{d \times r}$  are the optimal rank-$r$ orthonormal modes and their coefficients, respectively. The unresolved error is  the error of the optimal rank-$r$ approximation: $e_u(t)=\| \tilde{\bm{U}}(x,t) \tilde{\bm{Y}}(t)^T  - \bm{V}'(x,t) \|_F$, that is a direct result of truncating the $d-r$ least energetic modes. Thus, the optimal rank-$r$ approximation is obtained by minimizing:
\begin{equation}\label{eq:functional_opt}
\mathcal{E}_u(\tilde{\bm{U}}(x,t), \tilde{\bm{Y}}(t)) = \left\Vert  \tilde{\bm{U}}(x,t) \tilde{\bm{Y}}(t)^T  - \bm{V}'(x,t) \right\Vert_{F},
\end{equation}
subject to the orthonormality condition of $\tilde{\bm{U}}(x,t)$ modes. The optimal decomposition can be   obtained by performing instantaneous SVD of the sensitivity matrix, where  $\tilde{\bm{U}}(x,t)$ is the matrix of  $r$ most dominant left singular vectors of $\bm{V}'(x,t)$ and   $\tilde{\bm{Y}}(t) =  \tilde{\bm{Z}}(t)\tilde{\bm{\Sigma}}(t)$, where $ \tilde{\bm{Z}}(t) \in \mathbb{R}^{d \times r}$ and $ \tilde{\bm{\Sigma}}(t) =\mbox{diag}(\tilde{\bm{\sigma}}_1(t),\tilde{\bm{\sigma}_2}(t), \dots, \tilde{\bm{\sigma}}_r(t))$ are the matrix of the $r$ most dominant right singular vectors and the matrix of  singular values, respectively. It is straightforward to show that: $e_u(t) = (\sum_{i=r+1}^d \tilde{\sigma}_i^2(t))^{1/2}$.  The error $e_u(t)$ represents the minimum error that any rank-$r$ approximation can achieve, and therefore, it amounts to a lower bound for the f-OTD error: $e(t) \geq e_u(t)$. On the other hand, as with any reduced order model of a time-dependent system, the unresolved subspace induces a \emph{memory error} in the f-OTD approximation. This means that the unresolved error \emph{drives} the resolved error $e_r(t)$, and under appropriate conditions, it has been shown that for similar time-dependent basis low-rank approximations, $e_r(t)$ can be bounded  by: $e_r(t) \leq c_1 e^{c_2t}\int_{t_0}^t e_u(s)ds$  \cite{KL07} for $c_1,c_2>0$. The interplay between $e_u(t)$ and $e_r(t)$ can be more rigorously studied within the Mori-Zwanzig formalism \cite{CHK02}. These error estimates can guide an adaptive f-OTD, in which modes are added or removed to maintain the error below some threshold value \cite{Babaee:2017aa}, however these aspects are not in the scope of this paper and are not explored any further here. 
Since sensitivities can either be very small or very large with errors following the same trend, we  compute the relative error percentages as shown here:
\begin{equation}\label{eq:error}
    \text{\% Error} = \frac{e(t)}{\left\Vert\bm{V}'(x,t)\right\Vert_{F}} \times 100.
\end{equation}
Similar quantities are computed for $e_u(t)$ and $e_r(t)$. 
\subsection{Mode ranking}
In this section we present a procedure to rank the f-OTD modes and their coefficients according to their significance. To this end,  we start by considering the reduced correlation matrix $\bm{C}(t)$, which is in general a full matrix. This implies that the sensitivity coefficients are correlated and there exists a linear mapping from the correlated coefficients, $\bm{Y}(t)$, to the uncorrelated coefficients, $\hat{\bm{Y}}(t)\bs{\Sigma}(t)$, where $\hat{\bm{Y}}(t)$ are the orthonormal coefficients and $\bs{\Sigma}(t) = \mbox{diag}(\sigma_1(t),\sigma_2(t),\dots,\sigma_r(t))$ is a diagonal matrix of singular values. To find such a mapping, we consider the eigen-decomposition of $\bm{C}(t)$ as follows:
\begin{equation}
    \bm{C}(t)\bm{R}(t)=\bm{R}(t)\boldsymbol{\Lambda}(t),
\end{equation}
where $\bm{R}(t) \in \mathbb{R}^{r \times r}$ is a matrix whose columns contain the eigenvectors of $\bm{C}(t)$ and $\boldsymbol{\Lambda}(t)$ = diag($\lambda_1(t),\lambda_2(t),\dots,\lambda_r(t)$) is a diagonal matrix containing the eigenvalues of $\bm{C}(t)$. Since $\bm{C}(t)$ is a symmetric positive matrix, the matrix $\bm{R}(t)$ is an orthonormal matrix, i.e. $\bm{R}(t)^T \bm{R}(t) = \bm{I}$, and the eigenvalues are all non-negative and can be sorted as: $\lambda_1(t) > \lambda_2(t) > \dots > \lambda_r(t) \geq 0$. It is also straightforward to show that the singular values of the f-OTD low-rank approximation are $\sigma_i(t) = \lambda_i(t)^{1/2}$, for $i=1,2, \dots, r$. 
The ranked f-OTD components can be defined as:
\begin{equation*}
    \hat{\bm{Y}}(t) = \bm{Y}(t)\bm{R}(t)\bs{\Sigma}^{-1}(t), \quad \quad \hat{\bm{U}}(x,t) = \bm{U}(x,t)\bm{R}(t),
\end{equation*}
where the columns of $\hat{\bm{Y}}(t)$ and $\hat{\bm{U}}(x,t)$ are ranked by energy ($\sigma_i^2$) in descending order. We shall refer to \{$\hat{\bm{Y}}(t), \bs{\Sigma}(t),\hat{\bm{U}}(x,t)$\} as the bi-orthonormal form of the reduction. Since the above equations are simply an in-subspace rotation, \{$\hat{\bm{Y}}(t)\bs{\Sigma}(t),\hat{\bm{U}}(x,t)$\} and \{$\bm{Y}(t),\bm{U}(x,t)$\} yield equivalent low-rank approximations of the full-dimensional dynamics. This is easily verified by considering the bi-orthonormal form of the low-rank approximation as $\hat{\bm{U}}(x,t)\bs{\Sigma}(t)\hat{\bm{Y}}(t)^T$ $= \bm{U}(x,t)\bm{Y}(t)^T$, where we have made use of the identity $\bm{R}(t)^T \bm{R}(t) = \bm{I}$. We refer to $\hat{\bm{Y}}$ as the \emph{hidden} parametric space as each column of matrix $\hat{\bm{Y}}$ can be taken as a new ranked parameter that represents the contribution of all parameters ($\bs{\alpha}$).   In the following sections, all figures will be presented in bi-orthonormal form.

\section{Demonstration Cases}
\subsection{R\"{o}ssler system}
We first present a simple demonstration of  f-OTD  by computing sensitivities of the R\"{o}ssler system. The R\"{o}ssler system is governed by:
\begin{equation*}
    \frac{dv_1}{dt} = -v_2 - v_3, \quad \quad
    \frac{dv_2}{dt} = v_1 + \alpha_1v_2, \quad \quad
    \frac{dv_3}{dt} = \alpha_2 + v_3(v_1 - \alpha_3).
\end{equation*}
In the above equations, we set $\alpha_1=\alpha_2=0.1$ and $\alpha_3=14$, which are common values used to study the chaotic behavior of the attractor. The goal is to calculate the sensitivity of  $\bm{v}$ with respect to the model parameters $\bs{\alpha} = (\alpha_1, \alpha_2,\alpha_3)$ as $\partial\bm{v}/\partial\bs{\alpha}$. To this end, we take the derivative of the above system of equations with respect to model parameter $\alpha_i$ to obtain the sensitivity equation 
\begin{equation}\label{eqn:sensitivity_finite}
    \frac{d\bm{V}'}{dt} = \bm{L} \bm{V}' + \bm{F}',
\end{equation}
where
\[
\bm{L} = 
\begin{bmatrix}
    0 & -1 & -1 \\
    1 & \alpha_1 & 0 \\
    v_3 & 0 & v_1-\alpha_3
\end{bmatrix},
\quad
\bm{V}' = 
\begin{bmatrix}
    \vline & \vline & \vline \\
    \bm{v}'_1 & \bm{v}'_2 & \bm{v}'_3 \\
    \vline & \vline & \vline
\end{bmatrix},
\quad
\bm{F}' = 
\begin{bmatrix}
    0 & 0 & 0\\
    v_2 & 0 & 0 \\
    0 & 1 & -v_3
\end{bmatrix},
\]
and $\bm{v}'_i$ is the sensitivity of the position with respect to $\alpha_i$ and  $\bm{L} \in \mathbb{R}^{n\times n}$ and $\bm{F}' \in \mathbb{R}^{n\times d}$. 
\begin{figure}
\centering
\subfigure[]{
\includegraphics[width=.45\textwidth]{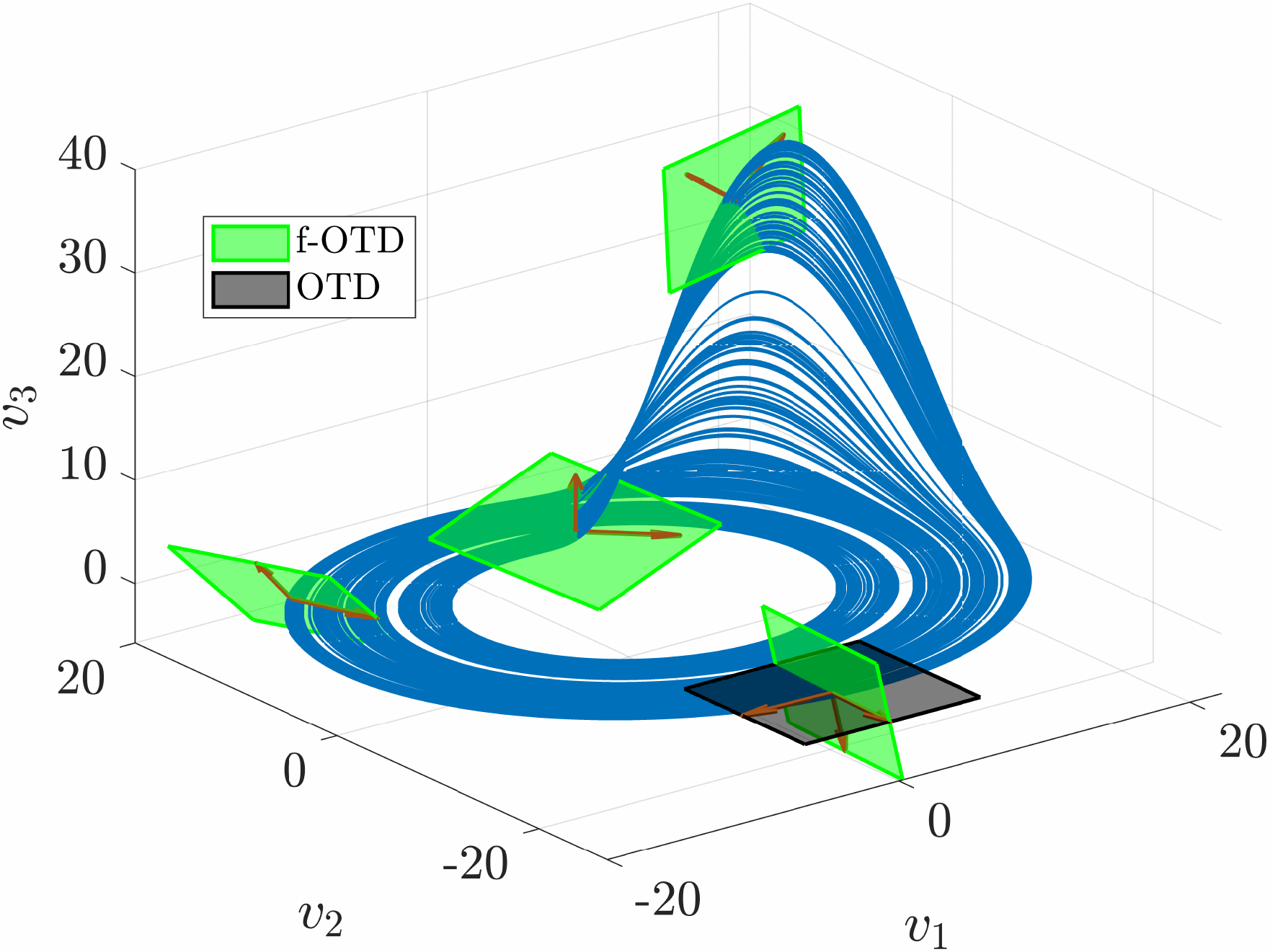}
\label{fig:ros_att}
}
\subfigure[]{
\includegraphics[width=.45\textwidth]{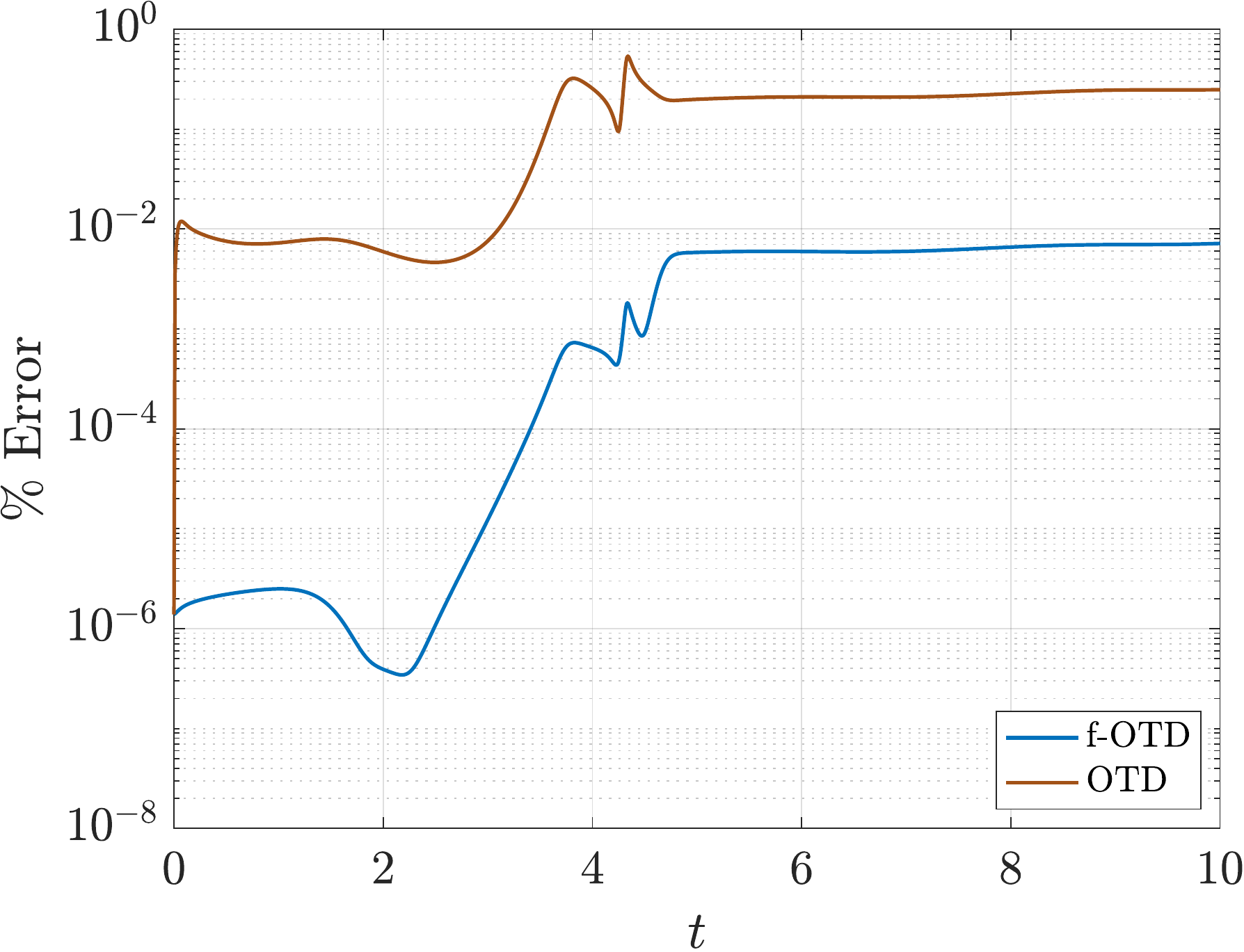}
\label{fig:ros_error}
}
\caption{(a) Chaotic R\"{o}ssler attractor with optimal f-OTD subspace shown in green and OTD subspace shown in black for $r=2$. Red arrows depict the orthonormal sensitivity vectors that define each subspace.  (b) Percent error for $e(t)$ plotted versus time for the f-OTD and OTD subspaces.}
\end{figure}
We choose a subspace with dimension $r=2$ for the low-rank approximation of the three-dimensional ($d=3$) sensitivities ($\bm{V}'$).  Although it is obvious that OTD modes are not based on parametric sensitivities and they are based on perturbations in the initial condition (IC) in all directions of the phase space, we believe it is instructive to contrast the OTD versus f-OTD to better understand f-OTD. To this end, we build two real-time ROMs using OTD modes and f-OTD modes. In the case of OTD, we solve the OTD evolution equation and we project the forced sensitivity Equation \ref{eqn:sensitivity_finite} onto the OTD modes, resulting in:
\begin{equation*}
    \frac{d\bm{U}_{otd}}{dt} = (\bm{I}-\bm{U}_{otd}\bm{U}_{otd}^T)\bm{L}\bm{U}_{otd} \quad  \mbox{and} \quad \frac{d\bm{Y}_{otd}}{dt}=\bm{Y}_{otd}\bm{U}^T_{otd}\bm{L}^T\bm{U}_{otd} + \bm{F}'^T\bm{U}_{otd}.
\end{equation*}
We also solved the f-OTD evolution Equations \ref{eqn:evol_modes} and \ref{eqn:evol_coeff} for the finite-dimensional system.  Both OTD and f-OTD modes are initialized with the same subspace and the evolution equations are solved for $T_f = 10$ units of time. These subspaces are initialized by first solving the full-dimensional sensitivity Equation \ref{eqn:sensitivity_finite} for one $\Delta t=10^{-2}$ and then computing the OTD  and f-OTD subspaces as the first two left singular vectors of $\bm{V}'(x,t=\Delta t)$. In Figure \ref{fig:ros_att}, both OTD and f-OTD subspaces are visualized along with the attractor of the R\"{o}ssler system. The OTD subspace is shown at only one instant for clarity and that point corresponds to the case where the nonlinear dynamics is in the $v_1-v_2$ plane. At this point, the OTD subspace is oriented such that it nearly coincides with the $v_1-v_2$ plane. This result is to be expected since the OTD subspace follows the sensitivities associated with the perturbations in the IC and we know that  the IC-perturbed solutions will  lie on the \emph{same} attractor. On the other hand, the f-OTD subspace is correctly oriented along the most sensitive subspace for perturbations in the model parameters, i.e. $\delta\bs{\alpha}=(\delta \alpha_1,\delta \alpha_2,\delta \alpha_3)$, which lead to perturbations in the attractor itself. That is, the perturbed solutions lie on \emph{different} attractors which can readily be seen as $\delta \bs{\alpha}$ results in nonzero $\delta v_3$, despite $v_3\simeq 0$. This results in the f-OTD subspace having a large out-of-plane component in the $v_3$ direction, which the OTD subspace fails to capture in Figure \ref{fig:ros_att}. In Figure \ref{fig:ros_error}, the percent errors of $e(t)$ are shown for OTD and f-OTD, which confirms that f-OTD  performs significantly better than OTD. This simple example demonstrates that the OTD basis is not optimal  and may be inaccurate for reduced order modeling of the forced sensitivity equation. 

\subsection{Chaotic Kuramoto Sivashinsky equation}
\begin{figure}
    \centering
    \includegraphics[width=\textwidth]{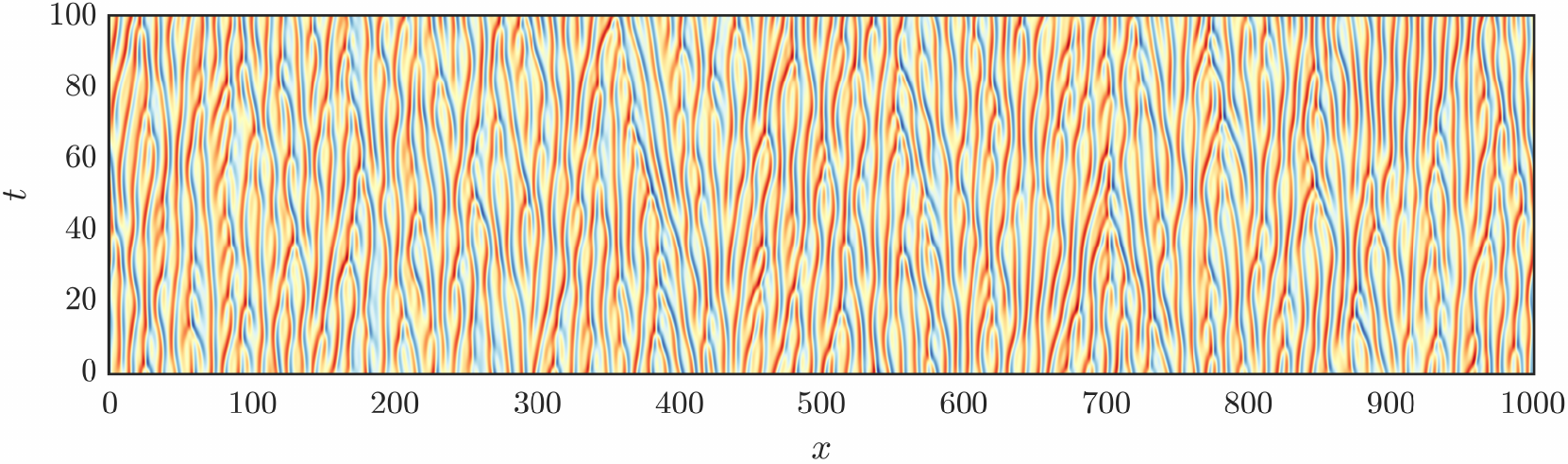}
    \caption{Solution of the chaotic Kuramoto-Sivashinksy equation $\bm{v}(x,t)$ solved on domain length $L=1000$ for $T_f=100$ units of time. }
    \label{fig:ks_nonlin}
\end{figure}

 The objective of this example is to evaluate the performance of f-OTD in computing sensitivities of a chaotic system with many  positive Lyapunov exponents and a high-dimensional parametric space. The intent of this example is not to compute the gradient of a time-averaged quantity for a chaotic system, but rather computing the solution of the sensitivity equation for a chaotic system with much larger unstable directions than the rank of the f-OTD subspace. For computing sensitivities of time-averaged quantities, one can use f-OTD in conjunction with Ruelle’s linear response formula \cite{R97,EHL04} to compute ensemble sensitivities.   To this end, 
we consider the sensitivity of the  Kuramoto Sivashinsky (KS) equation with respect to a time dependent forcing parameter $\alpha(t)$. The KS equation is a fourth order PDE given by:
\begin{align}\label{eqn:ks_pde}
    \pfrac{\bm{v}}{t} + \frac{1}{2}\pfrac{\bm{v}^2}{x} + \pfrac{^2 \bm{v}}{x^2} + \nu\pfrac{^4 \bm{v}}{x} = \alpha(t)\sin{\left( 2\pi x/L \right)}, \quad x\in[0,L],
\end{align}
where $\bm{v}=\bm{v}(x,t)$.  
Approximately 110 positive Lyapunov exponents exist for the parameters used in this study: $\nu=1$ and $L=1000$. 
The space time solution of Equation \ref{eqn:ks_pde} for these parameters is shown in Figure \ref{fig:ks_nonlin}. 
Here $\alpha(t)$ represents an infinite-dimensional parametric space. 


To compute the sensitivities numerically, we consider a discrete representation of $\alpha(t)$ in the interval $t_i\in[0,T_s]$, where $T_s\leq T_f$ is a subset of the full integration time $T_f$, and $t_i$ is a discrete instance in time. To this end,  we consider the value of $\alpha(t)$ at discrete time $t_i=(i-1)\times \Delta t$, where $\Delta t$ is the time step. This results in a vector, $\bs{\alpha} = (\alpha_1,\alpha_2,\dots,\alpha_d)$, where $\alpha_{i} = \alpha(t_i)$ and $d=T_s/\Delta t$ is the number of instances in time (i.e. number of parameters). In general, $\Delta t$ can be chosen independently of the numerical time integration step size, however, for simplicity, we use the same value of $\Delta t$ for both the parametric discretization and numerical integration of the nonlinear solver and f-OTD equations.   In this example, we consider $\Delta t = 10^{-2}$ and $T_s=10$, which results in $d=1000$ parameters. Further decrease in $\Delta t$ did not change our results. This leads to the sensitivity of $\bm{v}$ with respect to the value of $\alpha(t)$ at 1000 evenly spaced instances in time. We evolve these sensitivities over the interval $t\in[0,T_f]$ with $T_f=100$. We also choose $\alpha(t)=0$ for $t_i\in[0,T_f]$, and therefore, the nonlinear solver $\bm{v}(t)$ is the solution of the unforced KS equation. 
 
  
We consider the time-discrete form of Equation \ref{eqn:ks_pde} and differentiate with respect to design parameter $\alpha_i$. This leads to an evolution equation for the sensitivity of $\bm{v}$ with respect to $\alpha_i$, in which the linear operator and forcing terms are: 
\begin{equation}\label{eqn:ks_sens}
 \mathcal{L}(\bm{v}_i')= -\left[\pfrac{(\bm{v v}_i')}{x} + \pfrac{^2 \bm{v}'_i}{x^2} + \nu\pfrac{^4 \bm{v}'_i}{x}\right] \ \mbox{and} \ \bm{f}_i' = \delta(t-t_i)\sin\left( 2\pi x/L \right), \ i=1,2,\dots, d
\end{equation}
where $\delta(t-t_i)=0$ for $t\neq t_i$ and $\delta(t-t_i)=1$ for $t= t_i$. Our goal is to solve Equation \ref{eqn:ks_sens} using f-OTD.

\begin{figure}
\centering
\subfigure[]{
\includegraphics[width=.45\textwidth]{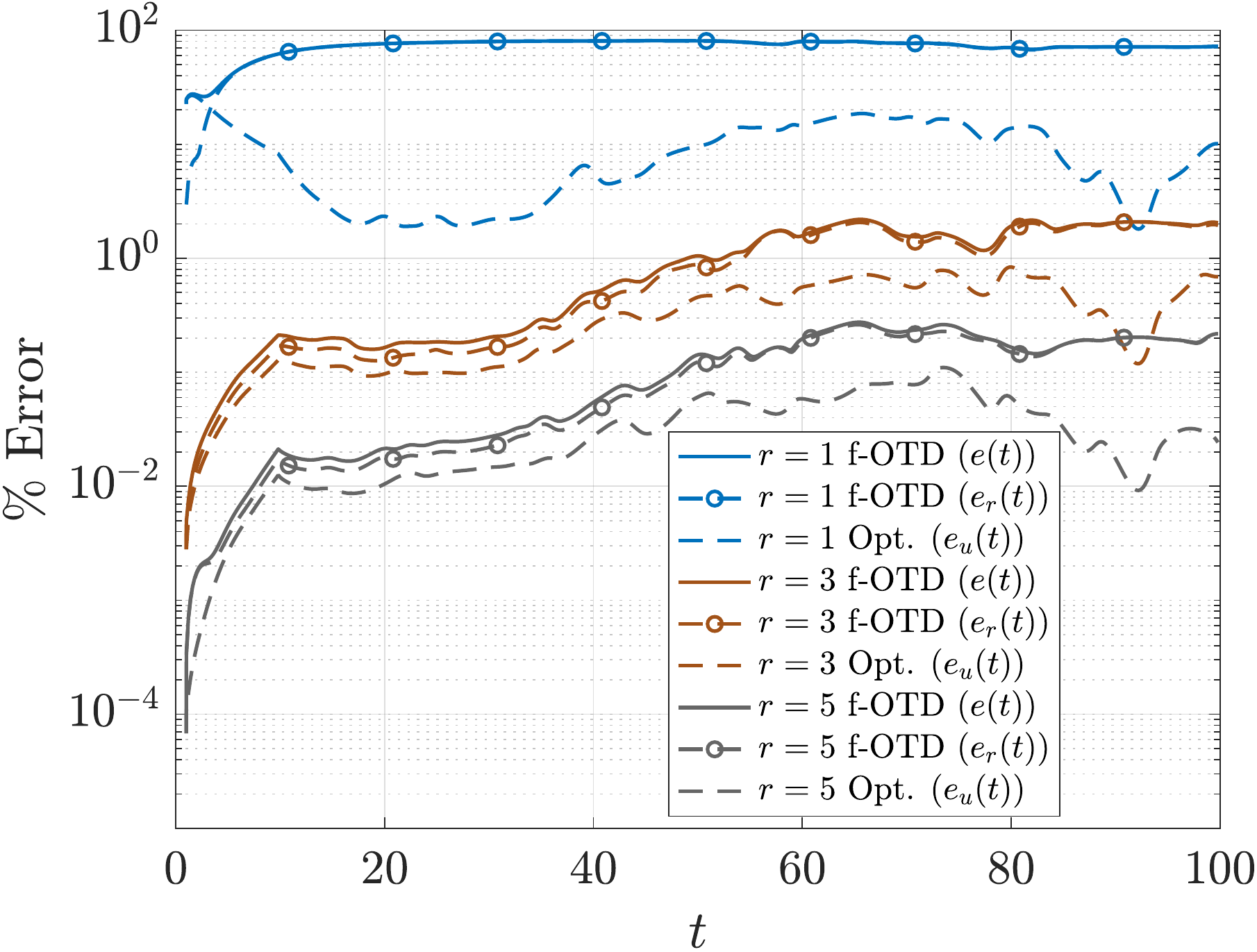}
\label{fig:ks_error}
}
\subfigure[]{
\includegraphics[width=.45\textwidth]{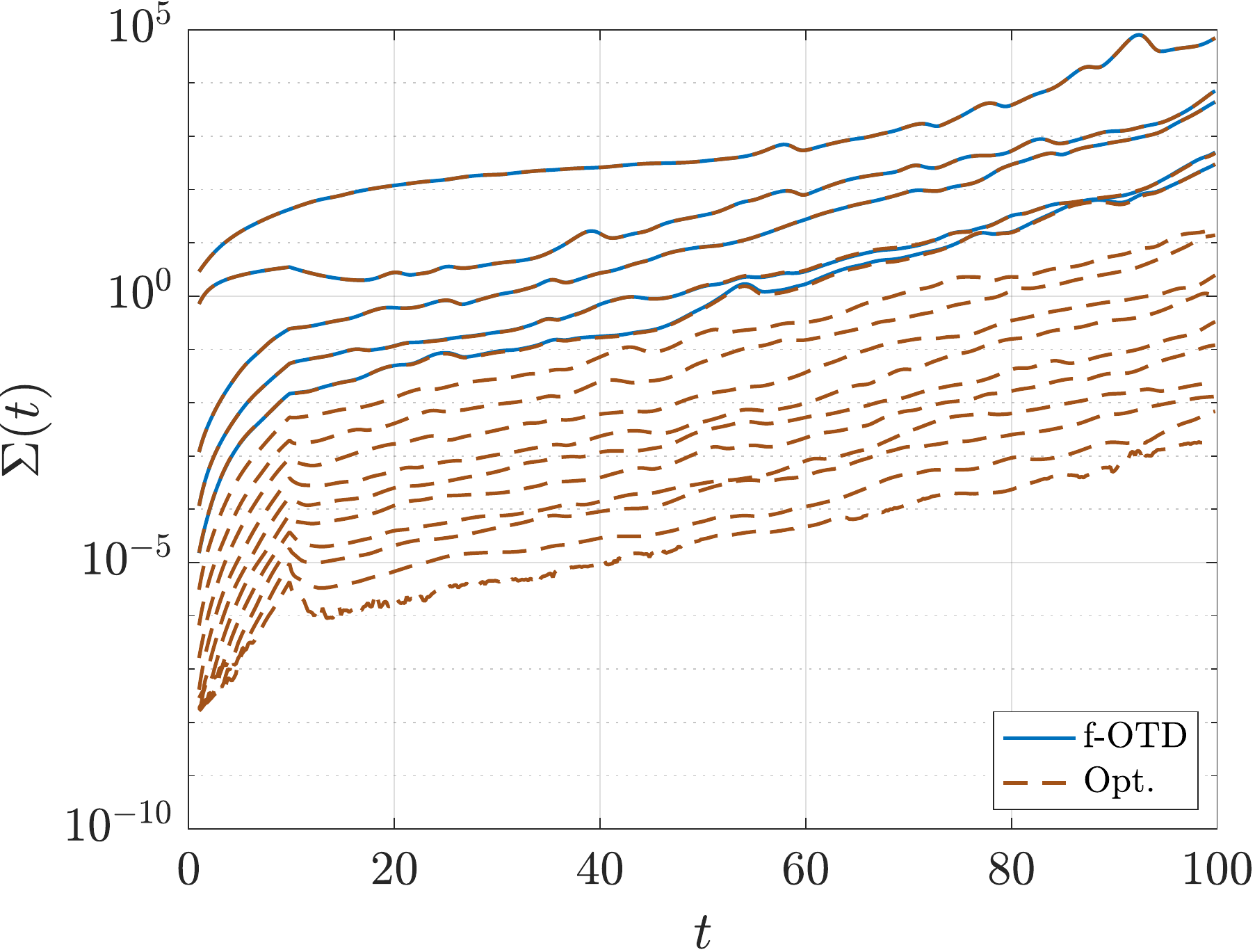}
\label{fig:ks_sing_val}
}
\caption{(a) Comparison of the reconstruction error between f-OTD approximation ($e(t)$) and optimal rank-$r$ approximation ($e_u(t)$) for different reduction sizes. Resolved error, $e_r(t)$, dominates the f-OTD error for long term integration. Error decreases as the number of modes increases. (b) Comparison of singular values between f-OTD and optimal low-rank decomposition for $r=5$.}
\end{figure}

We discretize the KS equation and the f-OTD equations using $n=2^{13}=8192$ Fourier modes and use exponential time-differencing Runge-Kutta fourth-order (ETDRK4) time stepping scheme \cite{KT05}. We verify our solution by directly solving Equation \ref{eqn:ks_sens} for all 1000 sensitivities. Further decreasing $\Delta t$ and increasing the number of Fourier modes did not change our results. We also compare the f-OTD error with that of optimal instantaneous same-rank approximation of the full sensitivities, which is obtained by computing the SVD of  $\bm{V}'(x,t)$ at each time. In Figure \ref{fig:ks_error}, we compare the reconstruction error of f-OTD ($e(t)$) with the reconstruction error of same-rank SVD ($e_u(t)$). We also show the resolved error $e_r(t)$, which measures the discrepancy between the f-OTD approximation and the optimal same-rank approximation. We compute these errors for $r=1,3$ and 5. While   the optimal low-rank approximation with a single mode captures approximately 99\% of the system energy of the full sensitivity (see Figure \ref{fig:ks_sing_val}),    the f-OTD approximation performs poorly with only a single mode, i.e., a dramatic reduction for 1000 sensitivities. This is a direct result of the memory effect from the lost interactions with the unresolved modes ($e_r(t)$) that ultimately dominate the error for long term integration.  By increasing the number of f-OTD modes, both $e(t)$ and $e_r(t)$ decrease. It is possible to control the error in real-time through an adaptive strategy that adds/removes modes with an appropriate criterion. For example, a candidate criterion could be  $p=\sigma^2_r(t)/\sum_{i=1}^r \sigma^2_i(t)$, where for $p<p_{th}$ the last mode can be removed and for $p>p_{th}$ a new mode can be added. See \cite{Babaee:2017aa} for similar strategies for adaptive mode addition and removal.   
In Figure \ref{fig:ks_sing_val}, we compare the 15 largest instantaneous singular values of quasimatrix $\bm{V}'(x,t)$ with those obtained from f-OTD with rank $r=5$, which shows  that f-OTD closely captures the most dominant subspace.

\begin{figure}[btp]
\centering
\subfigure[]{
\includegraphics[width=.45\textwidth]{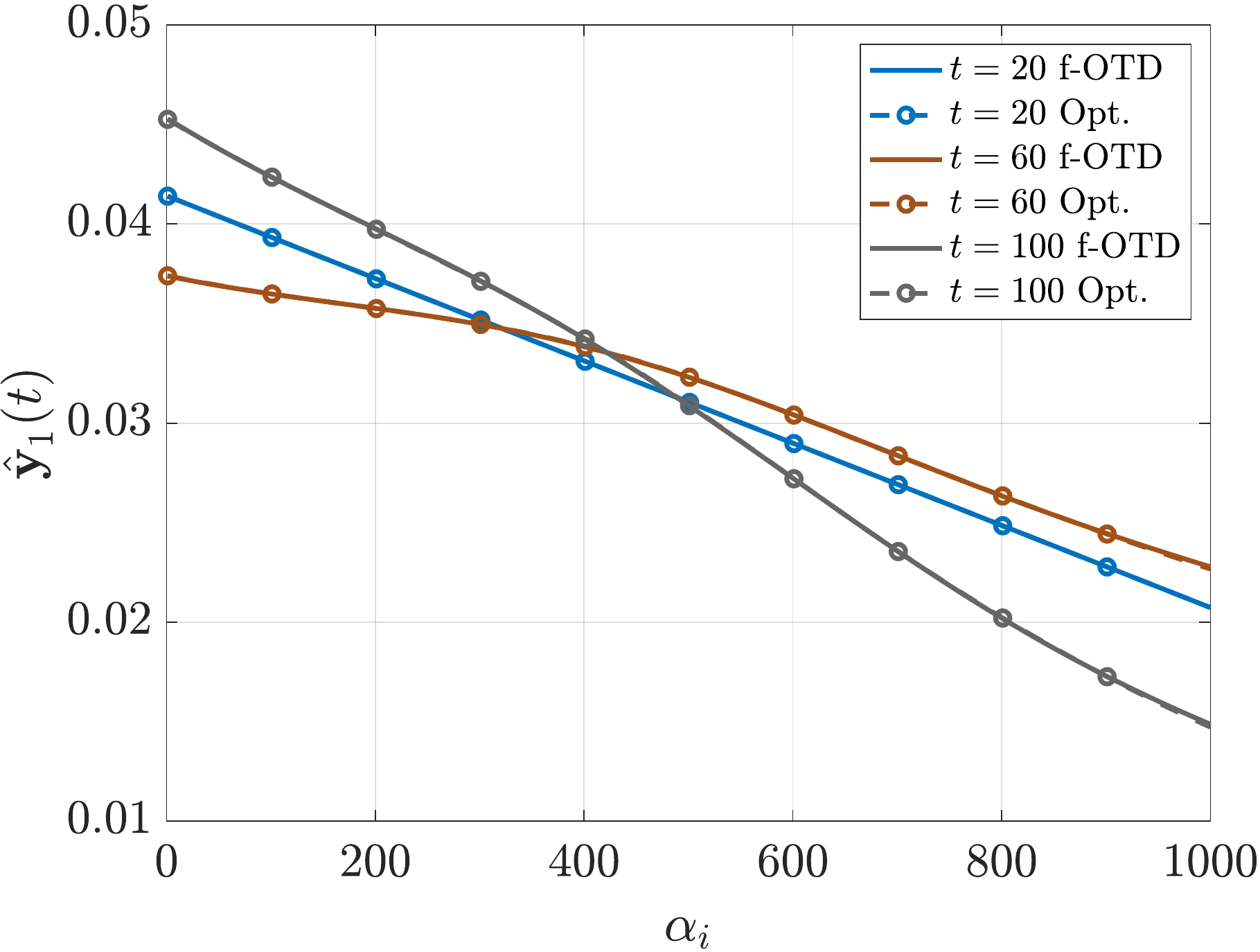}
\label{fig:ks_y1}
}
\subfigure[]{
\includegraphics[width=.45\textwidth]{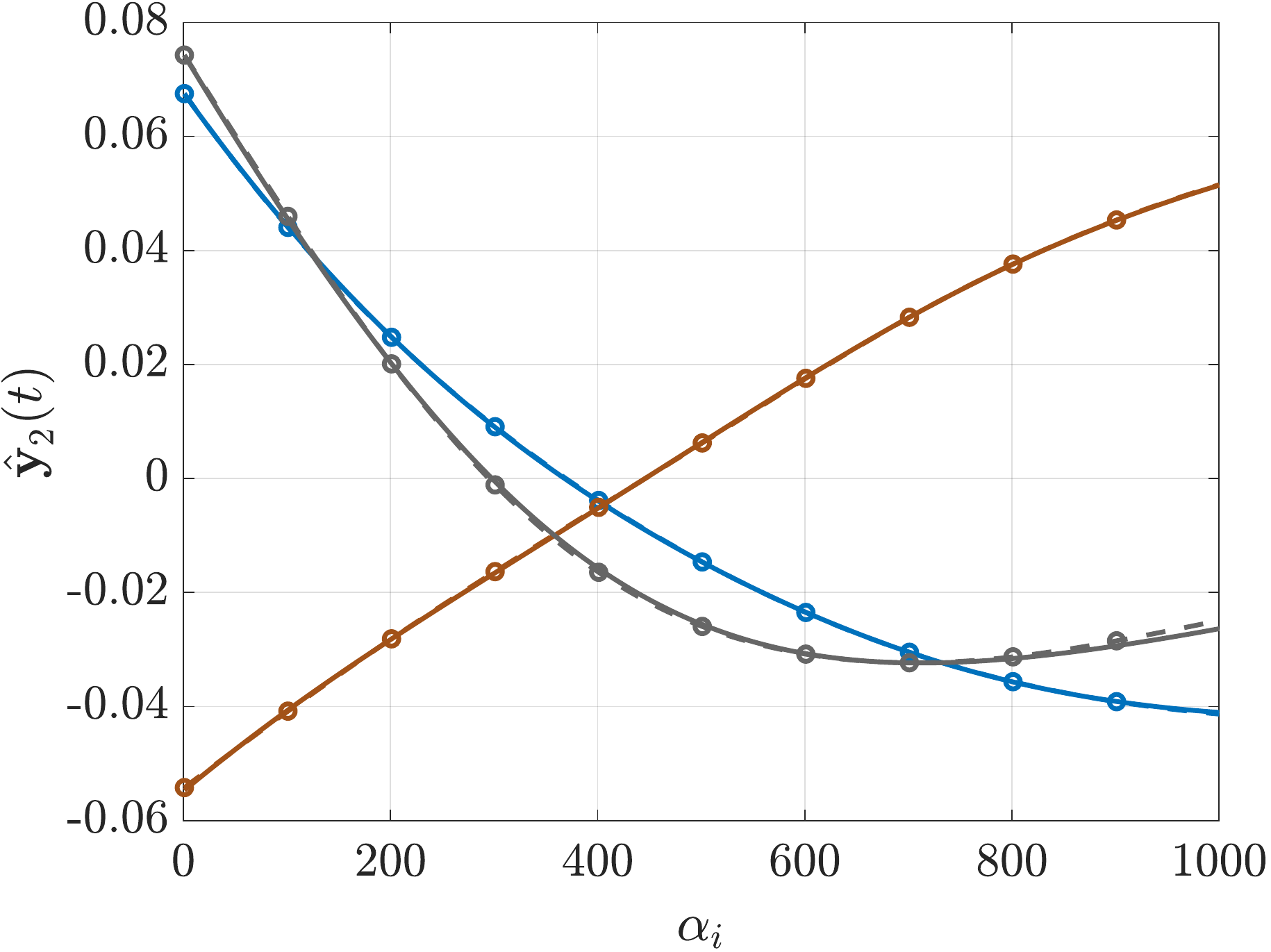}
\label{fig:ks_y2}
}
\caption{Kuramoto-Sivashinsky: The first two columns of the orthonormalized design variables matrix shown at different instances in time: (a) $\hat{\bm{y}}_1(t)$ (b) $\hat{\bm{y}}_2(t)$. The horizontal axis corresponds to the $i^{th}$ design parameter $\alpha_i$.}
\end{figure}

In Figures \ref{fig:ks_y1} and \ref{fig:ks_y2} the orthonormalized coefficients of the first two dominant f-OTD modes for the case of $r=5$ are compared to the right singular vectors from the instantaneous SVD of $\bm{V}'(x,t)$. These coefficients represent the hidden parametric space: for example, $\hat{\bm{y}}_1$ is a series of weights that represent the contribution of each of the $d=1000$ sensitivities to the most dominant direction of the full sensitivity matrix, $\hat{\bm{u}}_1$. Due to the chaotic nature of this problem, we observe that these coefficients can be highly time-dependent, especially for the lower energy modes; see $\hat{\bm{y}}_2$. Nevertheless, we have demonstrated that f-OTD extracts the most dominant subspace and associated coefficients of the sensitivity matrix for a chaotic system with large number of unstable directions and parameters.

\subsection{Species transport equation: Turbulent reactive flow}
In this example, we show how a single set of f-OTD modes can  lead to significant computational gains for computing sensitivities in problems with multiple coupled field variables,  where each field variable has a different linear  operator.  We consider a species transport problem, where parameter identification via sensitivity analysis  plays an important role in allocating computational and experimental resources to reduce parameter uncertainty.  Moreover, the sensitivity analysis is used to create reduced reaction mechanisms for complex chemical systems involving a large number of species and reactions. See references \cite{BOR15,LCR19,LLC20}.
\begin{figure}[btp]
    \centering
    \includegraphics[width=.7\textwidth]{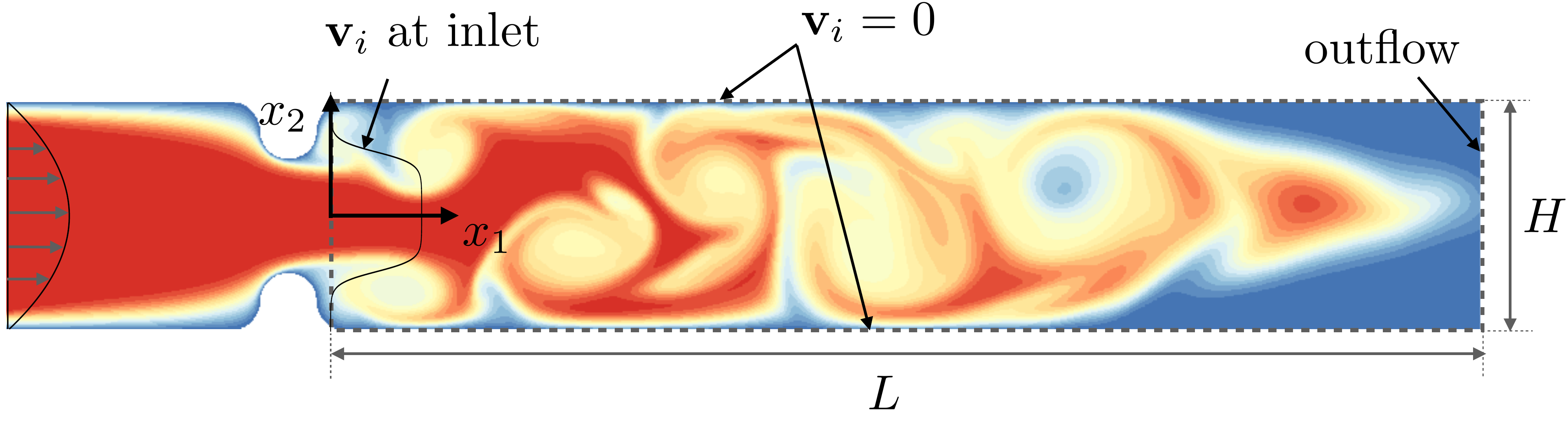}
    \caption{Schematic of the flow visualized with a passive scalar.}
    \label{fig:rxn_setup}
\end{figure}
\subsubsection{Problem setup}
To this end, we consider  a 2D incompressible turbulent reactive flow:
\begin{equation}\label{eq:ADR}
    \pfrac{\bm{v}_i}{t} + \left( \bm{w}\cdot\nabla \right)\bm{v}_i = \tilde{\kappa}_{ik}\nabla^{2}\bm{v}_{k} + \bm{s}_i,
\end{equation}
where $\bm{w}=(\bm{w}_{x_1}(x_1,x_2,t),\bm{w}_{x_2}(x_1,x_2,t))$ is the velocity field from the 2D incompressible Navier-Stokes equations, $\bm{v}_i=\bm{v}_i(x_1,x_2,t)$ is the concentration of species $i$, $\tilde{\kappa}_{ik}\in\mathbb{R}^{n_s\times n_s}$ is the diffusion coefficient matrix, and $\bm{s}_i=\bm{s}_i(\bm{v}_1,\bm{v}_2,\dots,\bm{v}_{n_s};\bs{\alpha})$ is the non-linear reactive source term. We choose a diagonal diffusion coefficient matrix, where the $i$th diagonal entry is the diffusion coefficient of the $i$th species, and $n_s$ is the number of species. For the reactive source term $\bm{s}_i$, we consider the  biological reactions used in \cite{LYTK15}. These terms are listed in Table \ref{tab:rxn} in Appendix \ref{app:D} for reference. A schematic of the flow is shown in Figure \ref{fig:rxn_setup}, where $L$ and $H$ are the channel length and height, respectively. The no-slip boundary condition is enforced at the top and bottom walls while the outflow boundary condition is enforced downstream. At the inlet a parabolic velocity with the average inlet velocity of  $\overline{w}$ is prescribed. The Reynolds number based on reference length of half the height ($H/2$) and the kinematic viscosity $\nu$  is $Re=\overline{w}H/2\nu=1000$.   The inlet boundary condition is $\bm{v}_i(0,x_2,t) = 1/2\big(\tanh{(x_2+H/2)/\delta} -\tanh{(x_2-H/2)/\delta}  \big)$ for all species, where $\delta = 0.1$.

The velocity field is governed by two-dimensional incompressible Navier-Stokes equation. We solved the velocity field once as it is independent from the species using spectral/hp elements method with 4008 quadrilateral elements and polynomial order 5. For more details on the spectral element method see for example \cite{KS05,Babaee:2013ab,Babaee:2013aa} . We then  solve the species transport equations and f-OTD equations in the rectangular domain shown by dashed lines in Figure \ref{fig:rxn_setup}. In the rectangular domain, we used structured spectral elements with 50 elements in $x_1$ direction and   15 elements in $x_2$ direction. We used spectral polynomial of order 5 in each direction. The velocity field was interpolated onto this grid.  The f-OTD equations, which are presented in the next sections,  and the species transport equation are integrated forward in time using RK4 with $\Delta t = 5 \times 10^{-4}$.

\subsubsection{f-OTD formulation}

\begin{figure}
    \centering
    \includegraphics[width=.8\textwidth]{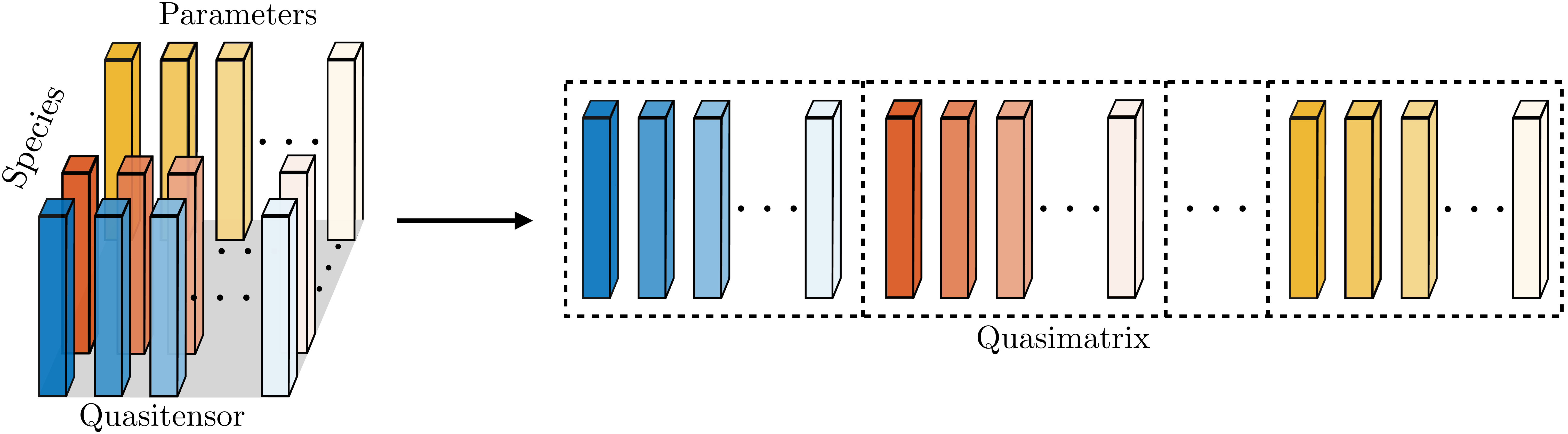}
    \caption{Schematic of the tensor flattening from a 3D quasitensor to a 2D quasimatrix.}
    \label{fig:tensor_flatten}
\end{figure}

Our goal is to calculate sensitivity of the species concentration with respect to the reaction parameters $\bs{\alpha}=(\alpha_1,\alpha_2,\dots,\alpha_{n_r})$, where $n_r$ is the number of reaction parameters. To this end, we take the derivative of the above equation with respect to reaction parameter $\alpha_j$ to obtain an evolution equation for the sensitivity:
\begin{equation}\label{eqn:rxn_sens_tensor}
    \pfrac{\tilde{\bm{v}}_{ij}'}{t} + \left( \bm{w}\cdot\nabla \right)\tilde{\bm{v}}_{ij}' =
    \tilde{\kappa}_{ik}\nabla^{2}\tilde{\bm{v}}_{kj}' + \tilde{\mathcal{L}}_{\bm{s}_{ik}} \tilde{\bm{v}}_{kj}'+\tilde{\bm{s}}'_{ij},
\end{equation}
where   $\tilde{\bm{v}}'_{ij}=\partial \bm{v}_i/\partial \alpha_j \in \mathbb{R}^{\infty\times 1}$  is the sensitivity of the concentration of species $\bm{v}_i$ with respect to reaction rate $\alpha_j$, $\tilde{\mathcal{L}}_{\bm{s}_{ik}}=\partial \bm{s}_i/ \partial \bm{v}_k$ is the linearized reactive source term, and $\tilde{\bm{s}}'_{ij}= \partial \bm{s}_i/ \partial \alpha_j$. In the above equation, $\tilde{\mathcal{L}}_{\bm{s}_{ik}} \tilde{\bm{v}}_{kj}'$ should be interpreted as a matrix-matrix multiplication for any  $(x_1,x_2)$ point in the physical space. In this notation,  sensitivities are represented by a quasitensor i.e. $\tilde{\bm{V}}' =[\bm{v}'_{ij}]$ with $i=1,2, \dots, n_s$ and $j=1,2,\dots, n_r$, where  $\tilde{\bm{V}}' \in \mathbb{R}^{\infty \times n_s \times n_r}$ is the third order quasitensor depicted in the left-hand side of Figure \ref{fig:tensor_flatten}.   Here $\tilde{\cdot}$ denotes terms associated with the tensor equation. In the discrete representation of   $\tilde{\bm{V}}'$, the dimension $\infty$ is replaced with the number of grid points. 

Solving for the sensitivities $\tilde{\bm{v}}_{ij}'$ using adjoint would require solving $n_s$ AEs:  one adjoint field for each species. See for example \cite{BOR15,LCR19,LLC20}.  
To solve for sensitivities using f-OTD, one could also solve for $n_s$ sets of f-OTD modes, i.e. one set of f-OTD modes for  each species. This straightforward approach would only exploit the correlation between sensitivities of  each species separately, i.e. correlations between $\bm{v}'_{ij}$ for a fixed $i$, while leaving  the correlations between sensitivities of different species  unexploited.    In this example, we demonstrate how a single set of f-OTD modes can be used to accurately model the entire sensitivity tensor. Therefore, the compression ratio both in terms of memory and computational cost in comparison to the full sensitivity equation is $r/d$. In comparison to AE, the compression ratio is $r/n_s$. Also, the f-OTD is a forward system and does not impose any I/O operation. To this end, we flatten the sensitivity tensor, as shown in Figure \ref{fig:tensor_flatten}, which results in a quasimatrix of size $\infty\times d$. Here, $d=n_s\times n_r$, where $n_s=23$ and $n_r=34$. This leads to a total of $d=782$ sensitivity equations that we seek to compute. In Appendix \ref{app:E}, we show that the flattened sensitivity evolution equation is:
\begin{equation}\label{eqn:rxn_sens_flatten}
    \pfrac{\bm{v}'_m}{t} + (\bm{w}\cdot \nabla)\bm{v}'_m = \kappa_{mn}\nabla^2 \bm{v}'_n + \mathcal{L}_{\bm{s}_{mn}}\bm{v}'_n + \bm{s}'_m,
\end{equation}
where $m(i,j)=j+(i-1)n_r$ and $n(i',j')=j' + (i'-1)n_r$, resulting in  $m,n=1,2,\dots,d$. Equation \ref{eqn:rxn_sens_tensor} is a tensor evolution equation, whereas Equation \ref{eqn:rxn_sens_flatten} is the equivalent matrix evolution equation.   The tensor flattening carried out here is similar to the unfolding carried out in the Tucker tensor decomposition \cite{KB09}. However, unlike Tucker tensor decomposition we do not consider flattening the tensor in the other two dimensions of species and parameters. Each $\bm{y}_k(t)$ is a vector of size $(n_s n_r)\times 1$ and contains coefficients for species and parameters.    Once the sensitivity tensor is flattened to a quasimatrix, we use f-OTD to extract low-rank structure from  the quasimatrix.  In Equation \ref{eqn:rxn_sens_flatten}, the linear operator changes from one species to the other due to the different diffusion coefficients $\kappa_{mn}$. In Appendix \ref{app:E} we show how f-OTD evolution equations can be derived for this case, which is different from the previous demonstration cases. 

\begin{figure}
\centering
\subfigure[]{
\includegraphics[width=.45\textwidth]{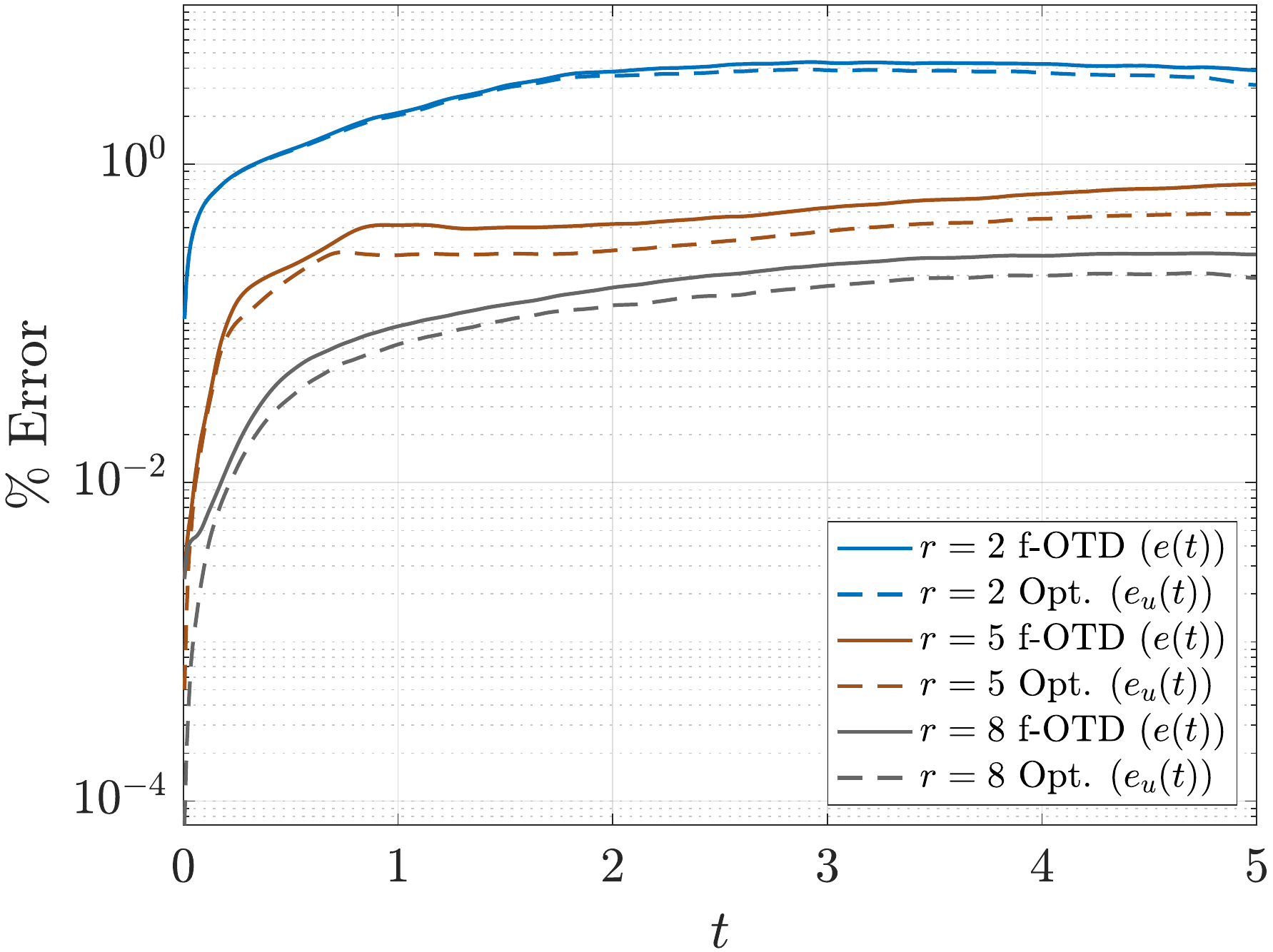}
\label{Fig:Jet_Error}
}
\subfigure[]{
\includegraphics[width=.485\textwidth]{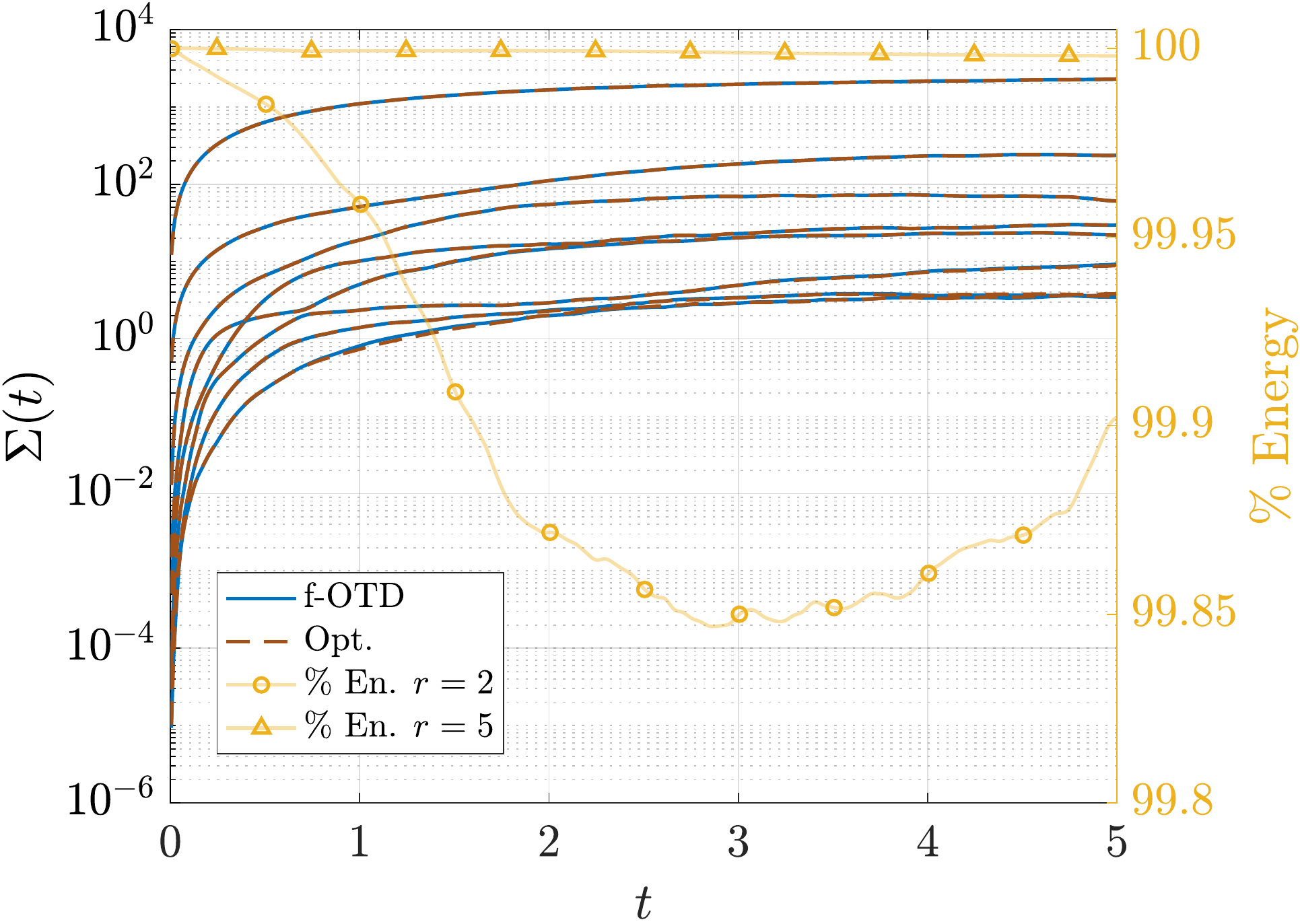}
\label{Fig:Jet_sing_vals}
}
\caption{(a) Percent error plotted as a function of time. Error decreases as the number of modes $r$ increases. (b) Singular values plotted as a function of time for $r=8$.}
\end{figure}

We solve Equations \ref{eq:uspec} and \ref{eq:yspec} for different f-OTD ranks along with the species transport (Equation \ref{eq:ADR}). In Figure \ref{Fig:Jet_Error} the  f-OTD error $(e(t))$ and optimal low-rank approximation error $(e_u(t))$ are shown using three different ranks of $r=2,5$ and 8. Again, we observe that the growth of $e(t)$ surpasses $e_u(t)$ for long term integration as a direct result of the lost interactions with the unresolved modes. However, with only 5-8 modes, we have shown that f-OTD can approximate 782 sensitivities with error on the order of 0.1\%. These results can be explained by studying Figure \ref{Fig:Jet_sing_vals}, where we observe that more than 99\% of the system energy is captured by the reduction. The \% energy is calculated from the singular values as \% En. $=\sum_{i=1}^{r}\sigma_i^2/\sum_{i=1}^{d}\sigma_i^2\times 100$, and can be used to get a sense of the dimensionality of the system, when expressed in the time-dependent basis. 
This allows the f-OTD algorithm to extract the latent features associated with the most dominant singular values and successfully approximate the full sensitivity tensor with a high degree of accuracy.
\begin{figure}[tbp]
    \centering
    \includegraphics[width=.9\textwidth]{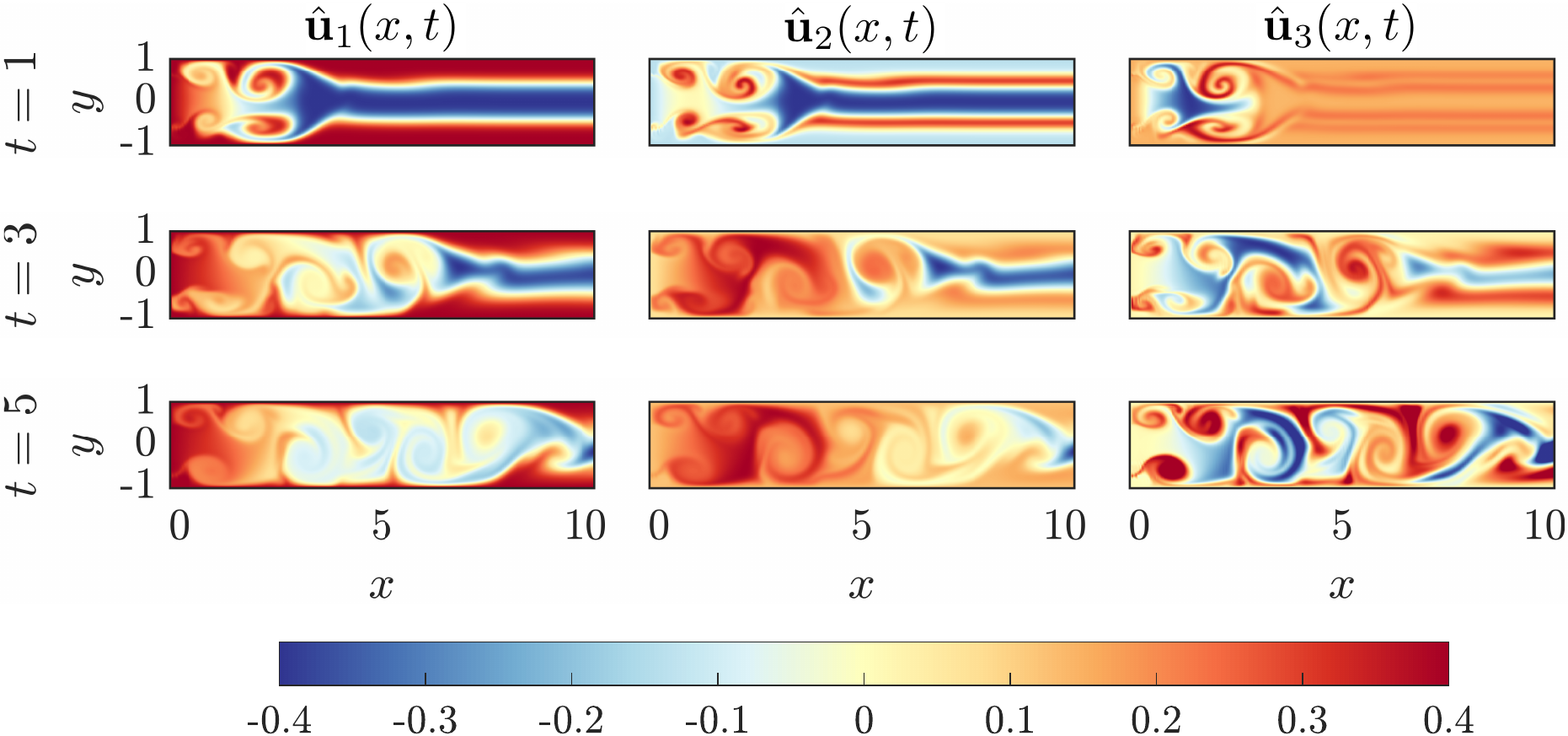}
    \caption{First three orthonormal f-OTD modes shown for $r=8$. Each row shows the modes at a different instance in time.}
    \label{fig:my_label}
\end{figure}

In Figure \ref{fig:my_label}, the time-dependent evolution of the three most dominant f-OTD modes are shown. These modes are energetically ranked where low mode numbers correspond to larger (higher energy) structures and high mode numbers correspond to finer (lower energy) structures in the flow. As opposed to static basis, such as POD or DMD,  the f-OTD modes evolve with the flow and exploit the instantaneous correlations between sensitivities. While this system is low-dimensional in the time-dependent basis, when expressed in POD or DMD basis, the system is high-dimensional and many modes are needed to capture the complex spatio-temporal evolution of $\bm{V}'$. See reference \cite{B19}  for comparison between time-dependent basis versus POD and DMD.  

\begin{figure}[hbt!]
    \centering
    \includegraphics[width=.9\textwidth]{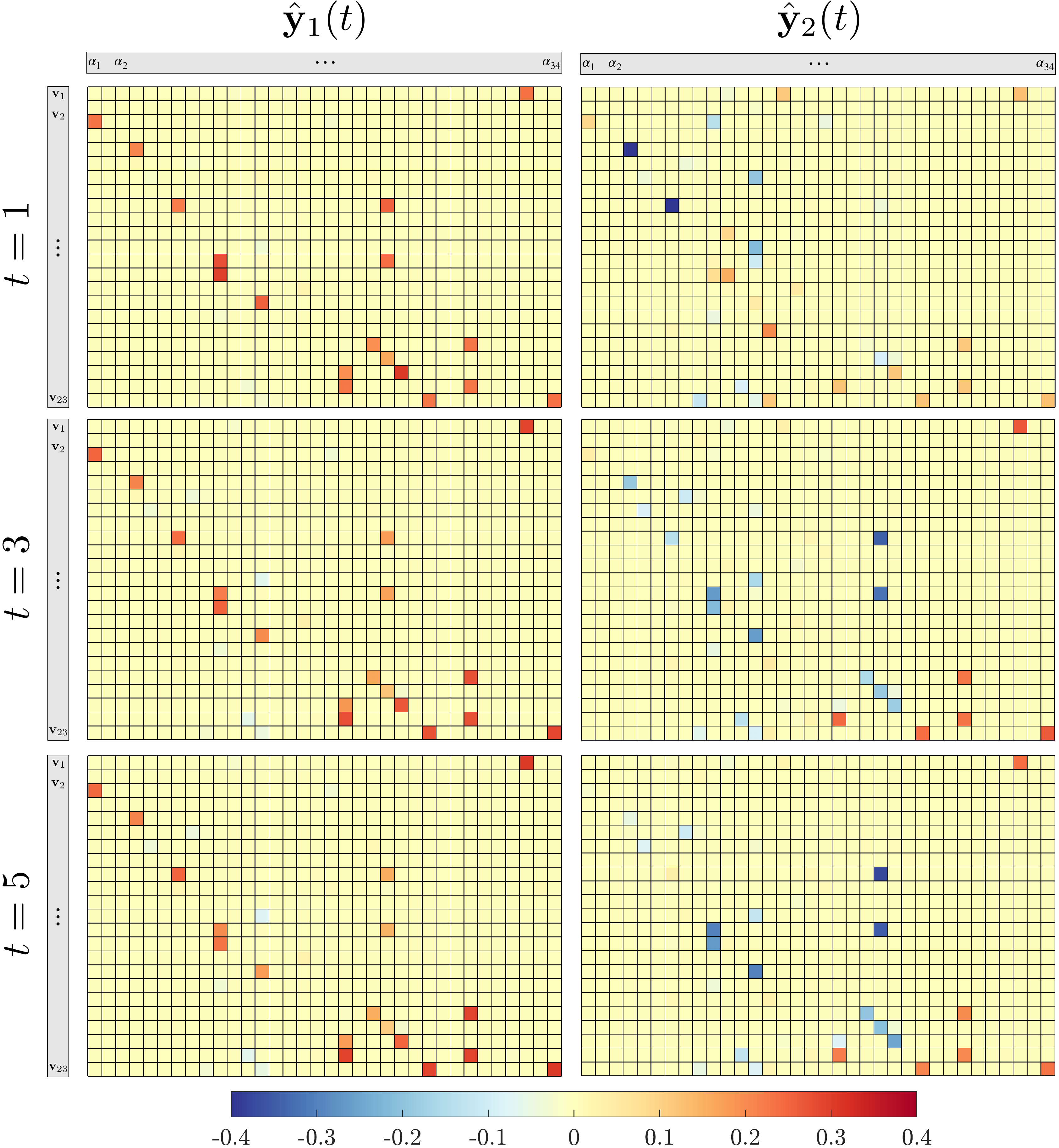}
    \caption{Orthonormalized f-OTD coefficients $\hat{\bm{y}}_1(t)$ and $\hat{\bm{y}}_2(t)$ visualized as a matrix with rows corresponding to species concentration and columns corresponding to reaction parameters. Color map shows most dominant sensitivities at different time instances.}
    \label{fig:y_matrix}
\end{figure}

To demonstrate the interpretability of the f-OTD decomposition,  we show how the hidden parameter space represented by $\hat{\bm{Y}}(t)$ can be used to identify the most important reaction parameters. In this context, importance refers to a parameter for which a small change in its value elicits a large change in the response of the system (i.e. highly sensitive).  To demonstrate this capability of f-OTD, the first two sensitivity coefficients are visualized as matrices in Figure \ref{fig:y_matrix}, where each $\hat{\bm{y}}_i$ is a $d\times 1$ vector that has been reshaped into an $n_s\times n_r$ matrix. In this form, each $\bm{v}'_{ij}$ is visualized using a heat map of the sensitivity coefficients, with rows corresponding to species $i$ and columns corresponding to reaction parameter $j$. Using this heat map, Figure \ref{fig:y_matrix} shows that only a handful of sensitivities are non-zero, while the majority have zero contribution for the entire duration of the simulation. 



\section{Conclusions}
We present a real time reduced order modeling approach for the computation of sensitivities in evolutionary systems governed by  time dependent ordinary/partial differential equations. The computational cost of solving the f-OTD equations of rank $r$ is roughly equivalent to that of solving $r$ forward sensitivity equations.  We demonstrated that the rank of f-OTD for two diverse applications is much smaller than the number of parameters. In contrast to adjoint based methods, f-OTD requires solving a system of  \emph{forward} equations and therefore it does not require any I/O operation. We showed that a single set of f-OTD modes can be formulated to compress  the sensitivities of multi-variable PDEs. We demonstrated this capability by computing sensitivities of multiple species with respect to reaction parameters in a turbulent reactive flow.  

 In contrast to traditional ROM approaches, f-OTD extracts the low-rank approximation  directly from the sensitivity equations as opposed to a data-driven approach, such as POD or DMD,  which requires the full-dimensional sensitivity data. Moreover, the low-rank subspace in the data-driven approach is fine tuned to particular operating conditions, whereas the f-OTD subspace is evolved with the dynamics of the system and does not require such fine-tuning. As such, f-OTD is  an \emph{on the fly} model compression that is achieved by extracting instantaneous correlated structures in the solution.


\appendix
\section{Optimality conditions of the variational principle}\label{app:A} 
For the sake of brevity, we forgo the explicit written dependencies on $x$, $t$, and $\bs{\alpha}$ in the following derivation. Using index notation, we start by expanding \ref{Eqn:min_princ}
\begin{align*}
    \mathcal{G}(\dot{\bm{U}},\dot{\bm{Y}}, \lambda) &= 
        \inner{\udot{i}}{\udot{j}}\left(\y{i}^T\y{j}\right) 
    +  \inner{\u{i}}{\u{j}}\left(\ydot{i}^T\ydot{j}\right) 
    +2 \inner{\udot{i}}{\u{j}}\left(\y{i}^T\ydot{j}\right) \\
    &-2 \inner{\udot{i}}{\mathcal{L}(\u{j})}\left(\y{i}^T\y{j}\right) 
    -2 \inner{\u{i}}{\mathcal{L}(\u{j})}\left(\ydot{i}^T\y{j}\right) \\
    &+  \inner{\mathcal{L}(\u{i})}{\mathcal{L}(\u{j})}\left(\y{i}^T\y{j}\right) 
    -2 \inner{\udot{i}}{\bm{F}'\y{i}} 
    -2 \inner{\u{i}}{\bm{F}'\ydot{i}} \\
    &+2 \inner{\mathcal{L}(\u{i})}{\bm{F}'\y{i}} 
    +  \big\|\bm{F}'\big\|_F^2 
    +  \lambda_{ij}\left(\inner{\u{i}}{\udot{j}}-\phi_{ij}\right).
\end{align*}
The first order optimality condition requires the derivative of $\mathcal{G}$ with respect to  $\dot{\bm{U}},\dot{\bm{Y}}$ and $ \lambda$ vanish. The derivative of $\mathcal{G}$ with respect to $\lambda$ produces the time derivative of the orthonormality constraint given by Equation \ref{eq:orthodot}. Provided that the  f-OTD modes  are orthonormal at $t=0$, the time integration of Equation \ref{eq:orthodot} reproduces the orthonormality condition of the f-OTD modes for  $t>0$: $\inner{\bm{u}_i}{\bm{u}_j}=\delta_{ij}$.  To take the derivative  of $\mathcal{G}$ with respect to  $\Tilde{\dot{\bm{u}}}_k$ we use  the Fr\'{e}chet differential as follows:
\begin{equation*}
\mathcal{G}'|_{\dot{\bm{U}}} \triangleq \lim_{\epsilon \rightarrow 0 } \frac{\mathcal{G}(\dot{\bm{U}}+\epsilon \dot{\bm{U}}',\dot{\bm{Y}},\lambda)-\mathcal{G}(\dot{\bm{U}},\dot{\bm{Y}},\lambda)}{\epsilon}.
\end{equation*}
Using the above definition we have:
\begin{align*}
    \mathcal{G}'|_{\udot{k}} &= 
    2 \inner{\dot{\bm{u}}'}{\udot{j}}\left(\y{k}^T\y{j}\right) 
    +2 \inner{\dot{\bm{u}}'}{\u{j}}\left(\y{k}^T\ydot{j}\right) 
    -2 \inner{\dot{\bm{u}}'}{\mathcal{L}(\u{j})}\left(\y{k}^T\y{j}\right) \\
    &-2 \inner{\dot{\bm{u}}'}{\bm{F}'\y{k}} 
    +  \lambda_{jk}\inner{\dot{\bm{u}}'}{\u{j}} = 0.
\end{align*}
The above equation can be written as $\inner{\dot{\bm{u}}'}{\nabla_{\dot{\bm{u}}_k} \mathcal{G}}$ and we observe that for any arbitrary direction $\dot{\bm{u}}'$, we must satisfy $\nabla_{\dot{\bm{u}}_k} \mathcal{G}=\bm{0}$. This leads to the following condition: 
\begin{equation}\label{eqn:udot_optimality}
     \nabla_{\dot{\bm{u}}_k} \mathcal{G} = 2 \udot{j}\left(\y{k}^T\y{j}\right) +2 \u{j}\left(\y{k}^T\ydot{j}\right) -2 \mathcal{L}(\u{j})\left(\y{k}^T\y{j}\right) 
    -2 \bm{F}'\y{k} +  \lambda_{jk}\u{j} = \bm{0}. 
\end{equation}
To eliminate $\lambda_{jk}$, we take the inner product of $\u{l}$ with \ref{eqn:udot_optimality} to obtain
\begin{align*}
    \left\langle\u{l}, \nabla_{\dot{\bm{u}}_k} \mathcal{G}\right\rangle &=
    2\phi_{lj}(\y{k}^T\y{j}) 
    +2 \delta_{lj}\left(\y{k}^T\ydot{j}\right) 
    -2 \inner{\u{l}}{\mathcal{L}(\u{j})}\left(\y{k}^T\y{j}\right) \\
    &-2 \inner{\u{l}}{\bm{F}'\y{k}} 
    + \lambda_{jk}\delta_{lj} = 0,
\end{align*}
where we have used $\inner{\u{l}}{\udot{j}}=\phi_{lj}$ and $\inner{\u{l}}{\u{j}}=\delta_{lj}$.
Rearranging for $\lambda_{lk}$ gives
\begin{equation*}
    \lambda_{lk} = 2\left[-\phi_{lj}(\y{k}^T\y{j}) - \left(\y{k}^T\ydot{l}\right) + \inner{\u{l}}{\mathcal{L}(\u{j})}\left(\y{k}^T\y{j}\right) + \inner{\u{l}}{\bm{F}'\y{k}}\right].
\end{equation*}
Dividing \ref{eqn:udot_optimality} by 2 and substituting $\lambda_{lk}$ gives
\begin{equation*}
\left[\udot{j} - \mathcal{L}(\u{j}) + \inner{\u{l}}{\mathcal{L}(\u{j}})\u{l} -\phi_{lj}\bm{u}_l \right]\left(\y{k}^T\y{j}\right) - \bm{F}'\y{k} + \inner{\u{l}}{\bm{F}'\y{k}}\u{l} = \bm{0}.
\end{equation*}
Rearranging the above equation for $\udot{j}$ we get
\begin{equation*}
    \boxed{\udot{j} = \mathcal{L}(\u{j}) - \inner{\u{l}}{\mathcal{L}(\u{j})}\u{l} + \left[\bm{F}'\y{k}-\inner{\u{l}}{\bm{F}'\y{k}}\u{l}\right]C_{kj}^{-1}+\phi_{lj}\bm{u}_l,}
\end{equation*}
where $C_{kj}=\y{k}^T\y{j}$.
Similarly, the first order optimality condition of $\mathcal{G}$ with respect to $\ydot{k}$ requires that
\begin{align}
    \pfrac{\mathcal{G}}{\ydot{k}} = \inner{\u{k}}{\u{j}}\ydot{j}
    + \inner{\udot{j}}{\u{k}}\y{j} 
    - \inner{\u{k}}{\mathcal{L}(\u{j})}\y{j} 
    - \inner{\bm{F}'}{\u{k}} 
    = \bm{0} \nonumber
\end{align}
Again, we use $\inner{\u{k}}{\u{j}}=\delta_{kj}$ and $\inner{\udot{j}}{\u{k}}=-\phi_{jk}$. Rearranging for $\ydot{k}$ gives
\begin{equation*}
    \boxed{\ydot{k} = \inner{\u{k}}{\mathcal{L}(\u{j})}\y{j} + \inner{\bm{F}'}{\u{k}} + \phi_{jk}\y{j}.}
\end{equation*}

\section{Equivalence of reductions}\label{app:B}
We prove the equivalence by using the evolution equation for the ${\bm{U},\bm{Y}}$ and using the matrix differential equation for the rotation matrix $\bm{R}$ and recovering the evolution equations for ${\tilde{\bm{U}},\tilde{\bm{Y}}}$. To this end, we substitute $\bm{U}=\tilde{\bm{U}}\bm{R}$ and $\bm{Y}=\tilde{\bm{Y}}\bm{R}$ into the quasimatrix form of Equations \ref{eqn:evol_modes} and \ref{eqn:evol_coeff}. The evolution equation for the orthonormal modes becomes:
\begin{align*}
\dot{\bm{U}} &= \dot{\tilde{\bm{U}}}\bm{R} + \tilde{\bm{U}}\dot{\bm{R}} \\
&= \mathcal{L}(\tilde{\bm{U}})\bm{R} - \tilde{\bm{U}}\bm{R}\inner{\tilde{\bm{U}}\bm{R}}{\mathcal{L}(\tilde{\bm{U}})\bm{R}} + [\bm{F}'\tilde{\bm{Y}}\bm{R} - \tilde{\bm{U}}\bm{R}\inner{\tilde{\bm{U}}\bm{R}}{\bm{F}'\tilde{\bm{Y}}\bm{R}}] + \tilde{\bm{U}}\bm{R}\bs{\Phi}.
\end{align*}
Substituting $\dot{\bm{R}}=\bm{R}\bs{\Phi}-\tilde{\bs{\Phi}}\bm{R}$ and  solving for $\dot{\tilde{\bm{U}}}$ yields 
\begin{align*}
    \dot{\tilde{\bm{U}}} &= \big[ \mathcal{L}(\tilde{\bm{U}})\bm{R} - \tilde{\bm{U}}\bm{R}\inner{\tilde{\bm{U}}\bm{R}}{\mathcal{L}(\tilde{\bm{U}})\bm{R}} + [\bm{F}'\tilde{\bm{Y}}\bm{R} - \tilde{\bm{U}}\bm{R}\inner{\tilde{\bm{U}}\bm{R}}{\bm{F}'\tilde{\bm{Y}}\bm{R}}] \\
    &+ \tilde{\bm{U}}\bm{R}\bs{\Phi} - \tilde{\bm{U}}[\bm{R}\bs{\Phi}-\tilde{\bs{\Phi}}\bm{R} \big]\bm{R}^T.
\end{align*}
Simplifying the above equation and using $\inner{\tilde{\bm{U}}\bm{R}}{\cdot}=\bm{R}^T\inner{\tilde{\bm{U}}}{\cdot}$ and $\bm{R}^{-1}=\bm{R}^T$, since $\bm{R}$ is an orthonormal matrix results in:
\begin{align*}
    \dot{\tilde{\bm{U}}} = \mathcal{L}(\tilde{\bm{U}}) - \tilde{\bm{U}}\inner{\tilde{\bm{U}}}{\mathcal{L}(\tilde{\bm{U}})} + [\bm{F}'\tilde{\bm{Y}} - \tilde{\bm{U}}\inner{\tilde{\bm{U}}}{\bm{F}'\tilde{\bm{Y}}}]\tilde{\bm{C}}^{-1} + \tilde{\bm{U}}\tilde{\bs{\Phi}},
\end{align*}
where $\tilde{\bm{C}}=\bm{R}\bm{C}\bm{R}^T$ and $\tilde{\bm{C}}^{-1}=\bm{R}\bm{C}^{-1}\bm{R}^T$, where $\bm{C}$ and $\tilde{\bm{C}}$ are similar matrices and thus have the same eigenvalues.  Following a similar procedure, the evolution equation for the coefficients becomes:
\begin{align*}
    \dot{\bm{Y}} &= \dot{\tilde{\bm{Y}}}\bm{R} + \tilde{\bm{Y}}\dot{\bm{R}} \\
    &= \tilde{\bm{Y}}\bm{R}\inner{\mathcal{L}(\tilde{\bm{U}})\bm{R}}{\tilde{\bm{U}}}\bm{R} + \inner{\bm{F}'}{\tilde{\bm{U}}}\bm{R} + \tilde{\bm{Y}}\bm{R}\bs{\Phi}.
\end{align*}
Substituting $\dot{\bm{R}}=\bm{R}\bs{\Phi}-\tilde{\bs{\Phi}}\bm{R}$ and  solving for $\dot{\tilde{\bm{Y}}}$ yields
\begin{align*}
    \dot{\tilde{\bm{Y}}} &=  \big[\tilde{\bm{Y}}\bm{R}\bm{R}^T\inner{\mathcal{L}(\tilde{\bm{U}})}{\tilde{\bm{U}}}\bm{R} + \inner{\bm{F}'}{\tilde{\bm{U}}}\bm{R} + \tilde{\bm{Y}}\bm{R}\bs{\Phi} - \tilde{\bm{Y}}[\bm{R}\bs{\Phi}-\tilde{\bs{\Phi}}\bm{R}] \big]\bm{R}^T \\
    &= \tilde{\bm{Y}}\inner{\mathcal{L}(\tilde{\bm{U}})}{\tilde{\bm{U}}} + \inner{\bm{F}'}{\tilde{\bm{U}}}  + \tilde{\bm{Y}}\tilde{\bs{\Phi}}.
\end{align*}
Thus, we have shown that the evolution of $\{\bm{U}(x,t),\bm{Y}(t)\}$ and $\{\tilde{\bm{U}}(x,t),\tilde{\bm{Y}}(t)\}$ according to Equations \ref{eqn:evol_modes} and \ref{eqn:evol_coeff} are equivalent.

\section{Exactness of f-OTD}\label{app:C}
Start by substituting $\bm{V}'(x,t) = \bm{U}(x,t)\bm{Y}(t)^T$ into the quasimatrix form of Equation \ref{eqn:sensitivity}:
\begin{align*}
    \dot{\bm{U}}\bm{Y}^T + \bm{U}\dot{\bm{Y}}^T = \mathcal{L}(\bm{U})\bm{Y}^T + \bm{F}'.
\end{align*}
Next we substitute $\dot{\bm{Y}}^T$ from Equation \ref{eqn:evol_coeff}
\begin{equation*}
    \dot{\bm{U}}\bm{Y}^T + \bm{U}\left( \bm{L}_r\bm{Y}^T + \inner{\bm{U}}{\bm{F}'} \right) = \mathcal{L}(\bm{U})\bm{Y}^T + \bm{F}',
\end{equation*}
where we have used $\bm{L}_r(t) = \inner{\bm{U}(x,t)}{\mathcal{L}(\bm{U}(x,t))}$. We multiply by $\bm{Y}$ from right and rearrange to get 
\begin{equation*}
 \dot{\bm{U}}\bm{C} = \mathcal{L}(\bm{U})\bm{C} - \bm{U}\bm{L}_r\bm{C} + \left( \bm{F}'\bm{Y} - \bm{U}\inner{\bm{U}}{\bm{F}'\bm{Y}} \right), 
\end{equation*}
where we have used $\bm{C} = \bm{Y}^T\bm{Y}$. Finally, we multiply by $\bm{C}^{-1}$ from right to get
\begin{equation*}
     \dot{\bm{U}} = \mathcal{L}(\bm{U}) - \bm{U}\bm{L}_r + \left( \bm{F}'\bm{Y} - \bm{U}\inner{\bm{U}}{\bm{F}'\bm{Y}} \right)\bm{C}^{-1},
\end{equation*}
which is the same as the evolution Equation \ref{eqn:evol_modes} for the orthonormal basis.  Here we have shown that the evolution of $\bm{V}'(x,t)$ under Equation \ref{eqn:sensitivity} is equivalent to the evolution of $\bm{U}(x,t)$ under Equation \ref{eqn:evol_modes}. That is, when $r=d$, $\bm{Y}(t)^T$ is a linear transformation that exactly maps the orthonormal subspace $\bm{U}(x,t)$ to $\bm{V}'(x,t)$.


\section{f-OTD Derivation for Tensor Sensitivities}\label{app:E} 
We start by considering the third order quasitensor $\tilde{\bm{V}}'=[\tilde{\bm{v}}'_{ij}] \in \mathbb{R}^{\infty \times n_s \times n_r}$ that we seek to flatten into a quasimatrix $\bm{V}'=[\bm{v}'_m]\in\mathbb{R}^{\infty\times d}$. For ease of reference, we rewrite the tensor evolution Equation \ref{eqn:rxn_sens_tensor} below:
\begin{equation*}
    \pfrac{\tilde{\bm{v}}_{ij}'}{t} + \left( \bm{u}\cdot\nabla \right)\tilde{\bm{v}}_{ij}' =
    \tilde{\kappa}_{ik}\nabla^{2}\tilde{\bm{v}}_{kj}' + \tilde{\mathcal{L}}_{\bm{s}_{ik}} \tilde{\bm{v}}_{kj}'+\tilde{\bm{s}}'_{ij},
\end{equation*}
where $i,k=1,2,\dots,n_s$ and $j=1,2,\dots,n_r$.
We define the indices $m(i,j)=j+(i-1)n_r$ and $n(i',j')=j' + (i'-1)n_r$, where $i'=1,2,\dots,n_s$ and $j'=1,2,\dots,n_r$. In the above equation, the terms $\tilde{\bm{v}}_{ij}'$ and $\tilde{\bm{s}}'_{ij}$ are flattened by replacing the index pair $ij$ with the single index $m$: $\bm{v}'_{m(i,j)} = \tilde{\bm{v}}_{ij}'$  and  $\bm{s}'_{m(i,j)} = \tilde{\bm{s}}_{ij}'$. Next, we define a new diffusion coefficient matrix $\kappa_{mn}\in\mathbb{R}^{d\times d}$ such that the $m$th diagonal entry is equal the diffusion coefficient of the $i$th species. That is, $\kappa_{mn}$ is independent of parameter index $j$ and remains constant across all sensitivities of a given species $i$. Finally, the linearized reactive source term is defined as $\mathcal{L}_{\bm{s}_{m(i,j)n(i',j')}} = \tilde{\mathcal{L}}_{\bm{s}_{ii'}} \delta_{j j'}$, where $\delta_{jj'}$ is the Kronecker delta and $n$ is a dummy index corresponding to $\bm{v}'_n$. From this definition, $\delta_{jj'}$ results in non-zero contribution to the summation over $n$ only for sensitivities with respect to parameter $j'=j$. Putting this all together, the above equation can be written as:  
\begin{equation}\label{eqn:}
    \pfrac{\bm{v}'_m}{t} + (\bm{w}\cdot \nabla)\bm{v}'_m = \kappa_{mn}\nabla^2 \bm{v}'_n+ \mathcal{L}_{\bm{s}_{mn}}\bm{v}'_n + \bm{s}'_m,
\end{equation}
where $\mathcal{L}_{\bm{s}_{mn}}\bm{v}'_n$ should be interpreted as a matrix-vector multiplication for any  $(x_1,x_2)$ point in the physical space. As a result of the parametric dependence of the linear operator, Equations \ref{eqn:evol_modes} and \ref{eqn:evol_coeff} do not hold for the tensor flattened equation. Therefore, we must derive new evolution equations for the f-OTD modes and coefficients for tensor flattened quantities. Substituting the approximation $\bm{v}'_m=\sum_{i=1}^{r} \bm{u}_iY_{mi}$ into the above equation, it is straightforward to show that the evolution equations for the f-OTD modes and coefficients are:
\begin{align}\label{eq:uspec}
    \dot{\bm{u}}_i = &-\left[(\bm{w}\cdot\nabla)\bm{u}_i - \bm{u}_j\inner{\bm{u}_j}{(\bm{w}\cdot\nabla)\bm{u}_i}\right] 
    + \left[ \nabla^2\bm{u}_k -  \bm{u}_j\inner{\bm{u}_j}{\nabla^2\bm{u}_k} \right]Y_{nk}\kappa_{mn}Y_{ml}C^{-1}_{il} \nonumber \\
    & + \left[ \mathcal{L}_{\bm{s}_{mn}}\bm{u}_k - \bm{u}_j\inner{\bm{u}_j}{\mathcal{L}_{\bm{s}_{mn}}\bm{u}_k} \right]Y_{nk}Y_{ml}C^{-1}_{il} 
     + \left[ \bm{s}'_m - \bm{u}_j\inner{\bm{u}_j}{\bm{s}'_m} \right]Y_{ml}C^{-1}_{il}, 
\end{align}
and
\begin{align}\label{eq:yspec}
    \dot{Y}_{mj} = &-\inner{\bm{u}_j}{(\bm{w}\cdot \nabla)\bm{u}_i}Y_{mi} + \inner{\bm{u}_j}{\nabla^2\bm{u}_i}Y_{ni}\kappa_{mn} \nonumber \\
    &+\inner{\bm{u}_j}{\mathcal{L}_{\bm{s}_{mn}}\bm{u}_i}Y_{ni} 
    +\inner{\bm{u}_j}{\bm{s}'_m}, 
\end{align}
where $\bm{Y}=[Y_{mi}]$ and the indices $m,n=1,2,\dots,d$ and $i,j,k=1,2,\dots,r$.
\clearpage

\section{Reactive source term specification}\label{app:D}
\begin{table}[hbtp]
\centering
\def\arraystretch{1.25}%
\begin{tabular}{l c} 
\hline
 $\bm{s}_1 = (\alpha_1 [13][2])/(\alpha_2 + [2]) -\alpha_3[1][15]$  \\ 
 $\bm{s}_2 = -(\alpha_1 [13][2]/(\alpha_2 + [2])$   \\ 
 $\bm{s}_3 = (\alpha_4 [9][4]/(\alpha_5 + [4]) - \alpha_6[3] - (\alpha_7[17][3])/(\alpha_8 +[3]) $ \\ 
 $\bm{s}_4 = (\alpha_4 [9][4])/(\alpha_5 + [4])$   \\ 
 $\bm{s}_5 = (\alpha_9 [9][6])/(\alpha_{10} + [6]) - \alpha_{11}[5]- (\alpha_{12}[17][5])/(\alpha_{13} +[5]) $ \\ 
 $\bm{s}_6 = -(\alpha_9 [9][6]/(\alpha_{10} + [6])$   \\
 $\bm{s}_7 = (\alpha_{14} [24][8])/(\alpha_{15} + [8]) - \alpha_{16}[7][15] -\alpha_{17}[16][7]$   \\
 $\bm{s}_8 = -(\alpha_{14} [24][8])/(\alpha_{15} + [8])$   \\ 
 $\bm{s}_9 = (\alpha_{18} [25][10])/(\alpha_{19} + [10]) -\alpha_{20}[9][15]$  \\ 
 $\bm{s}_{10} = -(\alpha_{18} [25][10])/(\alpha_{19} + [10])$   \\ 
 $\bm{s}_{11} = (\alpha_{21} [9][12])/(\alpha_{22} + [12]) - (\alpha_{23}[21][11])/(\alpha_{24} +[11]) $ \\ 
 $\bm{s}_{12} = -(\alpha_{21} [9][12])/(\alpha_{22} + [12])$   \\
 $\bm{s}_{13} = (\alpha_{25} [9][14])/(\alpha_{26} + [14]) - \alpha_{27}[13][15] -\alpha_{28}[13][19]$   \\
 $\bm{s}_{14} = -(\alpha_{25} [9][14])/(\alpha_{26} + [14])$   \\
 $\bm{s}_{15} = -(\alpha_3[1] + \alpha_{16}[7] + \alpha_{20}[9] +\alpha_{27}[13])[15]$   \\
 $\bm{s}_{16} = -\alpha_{17}[16][7]$   \\ 
 $\bm{s}_{17} = (\alpha_{29} [9][18])/(\alpha_{30} + [18]) -\alpha_{31}[17][19]$  \\  
 $\bm{s}_{18} = -(\alpha_{29} [9][18](\alpha_{30} + [18])$   \\
 $\bm{s}_{19} = -\alpha_{31}[17][19] -\alpha_{28}[13][19]$   \\ 
 $\bm{s}_{20} =0$\\
 $\bm{s}_{21} = (\alpha_{32} [20][22])/(\alpha_{33} + [22]) -\alpha_{34}[21][23]$  \\ 
 $\bm{s}_{22} = -(\alpha_{32} [20][22])/(\alpha_{33} + [22])$   \\ 
 $\bm{s}_{23} = -\alpha_{34}[21][23]$ \\
 \hline
\end{tabular}
\caption{Reactive source terms with species concentration denoted by $[\cdot]$. Each $\bm{s}_i$ is scaled by $10^2$ for time scale adjustment with the flow and the parameter values are assigned as follows: $\alpha_1 = 2.54\times 10^{-2}$, $\alpha_2 = 160$, $\alpha_3 = 3.74\times10^{-5} $, $\alpha_4 = 0.449$, $\alpha_5 = 1.12 \times 10^{5}$, $\alpha_6 = 5.13 \times 10^{-4}$, $\alpha_7 = 2.36 \times 10^{-2}$, $\alpha_8 = 14.6$, $\alpha_9 = 6.24\times 10^{-2}$, $\alpha_{10} = 140.5$, $\alpha_{11} = 3.93 \times 10^{-4}$, $\alpha_{12} = 2.36 \times 10^{-2}$, $\alpha_{13} = 14.6$, $\alpha_{14} = 5.523$, $\alpha_{15} = 160$, $\alpha_{16} = 8.01 \times 10^{-4}$, $\alpha_{17} = 1.11 \times 10^{-3}$, $\alpha_{18} = 3.105$, $\alpha_{19} = 1060$, $\alpha_{20} = 1.65 \times 10^{-3}$, $\alpha_{21} = 8.177$, $\alpha_{22} = 3160$, $\alpha_{23} = 3.456$, $\alpha_{24} = 2.50 \times 10^{5}$, $\alpha_{25} = 1.80 \times 10^{-5}$, $\alpha_{26} = 50$, $\alpha_{27} = 3.70\times 10^{-6}$, $\alpha_{28} = 3.00 \times 10^{-8}$, $\alpha_{29} = 9.01 \times 10^{-2}$, $\alpha_{30} = 3190$, $\alpha_{31} = 1.52 \times 10^{-9}$, $\alpha_{32} = 2.77 \times 10^{-2}$, $\alpha_{33} = 18$, and $\alpha_{34} = 2.22 \times 10^{-4}$}
\label{tab:rxn}
\end{table}

\section*{Acknowledgments}
We gratefully acknowledge the support received from the NASA Transformational Tools and Technologies Project, grant no. 80NSSC18M0150. We would also like to thank Dr. Joseph Derlaga for his support and intellectual insight throughout the process.   

\bibliographystyle{unsrt}
\bibliography{HB,MD}
\end{document}